\documentclass[10pt,a4paper]{amsart} 

\usepackage[margin=0.85in]{geometry} 
\usepackage{amsmath,amssymb,amsthm,mathtools}
\usepackage{booktabs}
\usepackage[dvipsnames]{xcolor}
\usepackage{tikz}
\usepackage{pgfplots}
\usepackage{hyperref}

\pgfplotsset{compat=1.18}

\definecolor{SteinBlue}{RGB}{0,114,178}
\definecolor{MYZOrange}{RGB}{230,159,0}
\definecolor{GHLWGreen}{RGB}{0,130,105}
\definecolor{BoundaryPurple}{RGB}{145,70,135}
\definecolor{KnownFill}{RGB}{225,241,234}
\definecolor{GainFill}{RGB}{220,235,247}
\definecolor{GapFill}{RGB}{255,239,196}
\definecolor{ExcludedFill}{RGB}{249,224,224}

\newcommand{\R}{\mathbb R}
\newcommand{\M}{\mathcal M}
\newcommand{\eps}{\varepsilon}

\newtheorem{theorem}{Theorem}[section]
\newtheorem{proposition}[theorem]{Proposition}
\newtheorem{conjecture}[theorem]{Conjecture}
\newtheorem{corollary}[theorem]{Corollary}
\newtheorem{remark}[theorem]{Remark}
\newtheorem{lemma}[theorem]{Lemma}

\numberwithin{equation}{section}

\title[Stein's maximal spherical means]
{Stein's generalized spherical maximal operator}

\author[L.\ Roncal]{Luz Roncal}
\address[Luz Roncal]{
	BCAM -- Basque Center for Applied Mathematics,
	48009 Bilbao, Spain,
	\newline \phantom{\quad} \&
	Ikerbasque, Basque Foundation for Science,
	48011 Bilbao, Spain,
	\newline \phantom{\quad} \&
	Universidad del Pa\'is Vasco / Euskal Herriko Unibertsitatea,
	48080 Bilbao, Spain}
\email{\href{mailto:lroncal@bcamath.org}{lroncal@bcamath.org}}

\date{\today}

\begin{document}

\begin{abstract}
This survey concerns Stein's analytic family of generalized spherical means
and the associated maximal operator.  We review Stein's original theorem, the
sharp necessary conditions for boundedness, and the complete theory for radial
functions.  We then explain the connection between positive results and local
smoothing estimates for the half-wave propagator and, more generally, for
Fourier integral operators.  This point of view yields the sharp
two-dimensional result through the work of Guth, H. Wang, and R. Zhang and, when
combined with the recent local smoothing theorems, an improved
higher-dimensional range.  As a byproduct, the maximal operator
associated with the first radial derivative of spherical means is bounded
on $L^p(\R^4)$ for $2<p<3$, a case left open in recent work of
Kwa\'snicki. We keep track of the  relation between the
local smoothing exponents and the resulting maximal estimates, and identify
the region that remains open. 
\end{abstract}

\maketitle

\section{Purpose and scope of the survey}
\label{sec:intro}

The central object of this survey is Stein's \emph{generalized spherical
maximal operator}.  Fix $n\ge2$.  We use the normalization for which the
parameter $\alpha=0$ corresponds to normalized spherical surface measure.  For
$x\in\R^n$ and $\operatorname{Re}\alpha>0$, define
$$
 m_\alpha(x)=
 \frac{\Gamma(\alpha+n/2)}{\pi^{n/2}\Gamma(\alpha)}
 (1-|x|^2)_+^{\alpha-1},
 \qquad
 m_{\alpha,t}(x)=t^{-n}m_\alpha(x/t),
$$
and, initially for $f\in\mathcal S(\R^n)$, set
$$
 \M_t^\alpha f(x)=f*m_{\alpha,t}(x),
 \qquad
 M^\alpha f(x)=\sup_{t>0}|\M_t^\alpha f(x)|.
$$
The normalized family is meromorphic in $\alpha$ and is holomorphic for
$
 \alpha\in\mathbb C\setminus
 \{-n/2,-n/2-1,-n/2-2,\ldots\}$.
Throughout the survey, orders belonging to this exceptional set are excluded.
With the Fourier transform convention
$\widehat f(\xi)=\int_{\R^n}e^{-ix\cdot\xi}f(x)\,dx$, one has
\begin{equation}
\label{eq:multiplier}
 \widehat{m_\alpha}(\xi)
   =c_{n,\alpha}|\xi|^{1-n/2-\alpha}
     J_{n/2+\alpha-1}(|\xi|),
 \qquad
 c_{n,\alpha}=2^{\alpha+n/2-1}\Gamma(\alpha+n/2),
\end{equation}
where $J_\nu$ denotes the Bessel function of the first kind.  In particular,
$m_0$ is the normalized surface measure on $\mathbb S^{n-1}$.  As we shall
see, the oscillatory asymptotics of the Bessel function provide a bridge
between Stein's analytic family, the half-wave propagator, and local
smoothing.  The transfer principle in Section~\ref{sec:transfer} makes
precise how improvements in local smoothing yield larger $L^p$ boundedness
regions for $M^\alpha$.\footnote{Among the AI-generated
manuscripts made public by OpenAI on October 6, 2026
(\url{https://github.com/openai/math}), one claims the critical $L^3$ local
smoothing estimate, with arbitrarily small loss, for the wave equation in
three spatial dimensions, that is, the Euclidean local smoothing conjecture (see Subsection \ref{sub:conels}) for $n=3$; a
second manuscript, dated October 25, 2026 and accompanied by a Lean
formalization, claims the conjecture for every $n\ge2$.  At the time of
writing, neither claim has been peer reviewed, and the formal statements
have not been checked against the manuscripts by human reviewers.  This
survey is written accordingly: the Euclidean local smoothing conjecture is
treated as open.  Even if these claims are confirmed, several of the
questions discussed here, in particular the endpoint ones, remain; see
Section~\ref{sec:open}.}

The survey has three main goals.
\begin{enumerate}
\item[1)] Trace the development of the boundedness theory, beginning with Stein's
  paper \cite{Stein1976}.
\item[2)] Give a unified account of the transfer principle developed by Miao,
  Yang, and Zheng and later made explicit by Liu, Shen, Song, and Yan, with
  careful attention to its derivative losses.
\item[3)] Apply that principle to the strongest available local smoothing
  theorems, including the recent result by Gan, He, Li, and
Wu 
  \cite{GHLW2026} and its even-dimensional refinement by X.~Wang
  \cite{Wang2026}, and isolate the remaining gap.  A byproduct is
  Corollary~\ref{cor:n4}: the maximal operator of the first radial
  derivative of spherical means is bounded on $L^p(\R^4)$ for $2<p<3$,
  which settles the case $n=4$ of a question left open in
  \cite{Kwasnicki2026}.
\end{enumerate}

Our scope is deliberately narrower than that of the recent survey by Roos and
Seeger \cite{RoosSeeger2026}.  Here the focus is Stein's analytic family, and
the principal ingredient is the local smoothing theory of Fourier integral
operators (FIOs).

\section{The boundedness diagram}
\label{sec:diagram}

We begin by presenting the successive known boundedness regions for
$M^\alpha$ in the $(1/p,\operatorname{Re}\alpha)$-plane.  Placing these
regions in a common diagram gives an overview of the development of the
problem and fixes the principal benchmarks for the rest of the survey.

\subsection{Classical positive region}

In his seminal note \cite{Stein1976}, Stein introduced the analytic family
$M^\alpha$ as an auxiliary tool in the study of the classical spherical
maximal function.\footnote{To the best of our knowledge, the terminology
\emph{Stein's maximal spherical means} for $M^\alpha$ first appeared in
\cite{MYZ2017}.}  He proved the following theorem.

\begin{theorem}[\cite{Stein1976}]
\label{thm:Stein}
Let $n\ge3$.  Then $M^\alpha$ is bounded on $L^p(\R^n)$ provided that
\begin{equation}
\label{eq:Stein-region}
 \operatorname{Re}\alpha>1-n+\frac np,
 \quad 1<p\le2,
 \qquad\text{or}\qquad
 \operatorname{Re}\alpha>\frac{2-n}{p},
 \quad 2\le p\le\infty.
\end{equation}
\end{theorem}

The first condition in \eqref{eq:Stein-region}, corresponding to
$1<p\le2$, is optimal.  The planar case of the classical spherical maximal
operator, obtained by setting $\alpha=0$, was subsequently established by
Bourgain \cite{Bourgain1986}.  In his monograph, Stein observed that the
optimal range for $p>2$ and $n\ge2$ was ``still a mystery''
\cite[Ch.~XI, \S4.10]{SteinHA}.

\subsection{The radial theory and the conjectural boundary}

For real orders, the radial theory identifies the conjectural boundary and
distinguishes the behavior of its endpoints.  It was initiated by
Ciaurri, Nowak, and the author \cite{CNR2017,CNR2020} for the
two-parameter family of generalized spherical mean Radon transforms, in
which the dimension enters as a continuous parameter and the second
parameter is the order $\alpha$ of this survey.  On radial functions,
\cite{CNR2017} gives two-weight mixed norm estimates for the means, and
\cite{CNR2020} power-weighted estimates for the maximal operator, together
with almost everywhere convergence results for the Euler--Poisson--Darboux
and wave equations; a first version of Conjecture~\ref{conj:NRS-general}
was formulated there.  Set
\begin{equation}
\label{eq:qn}
 q_n:=\frac{2n}{n-1}.
\end{equation}
The following result was proved by Nowak, the author, and
Szarek.\footnote{They use $\beta$ for the order of the generalized spherical
mean.  We state their result using $\alpha$ in order to keep the notation
uniform throughout the survey.}

\begin{theorem}[{\cite[Corollary~1.5]{NRS2023}}]
\label{thm:NRS-radial}
Let $n\ge2$, let $\alpha\in\R$ satisfy
$\alpha>(1-n)/2$, and let $1\le p<\infty$.  Then $M^\alpha$ is bounded on
$L^p_{\mathrm{rad}}(\R^n)$ if and only if $p>1$ and
\begin{equation}
\label{eq:radial-region}
\begin{cases}
 \alpha>1-n+\dfrac np, &1<p<2,\\[5pt]
 \alpha>\dfrac1p-\dfrac{n-1}{2},
     &2\le p\le q_n,\\[5pt]
 \alpha\ge-\dfrac{n-1}{p},&p>q_n.
\end{cases}
\end{equation}
\end{theorem}

Any $L^p$ estimate for general functions restricts, of course, to radial
functions.  Consequently, Theorem~\ref{thm:NRS-radial} already gives the
necessary boundary for real parameters satisfying $\alpha>(1-n)/2$.  Notice
that this restriction on $\alpha$ does not remove any point from the predicted
positive region when $p\ge2$, since each of the lower bounds in
\eqref{eq:radial-region} lies strictly above $(1-n)/2$.  An earlier
contribution to the radial theory is due to Seeger, Wainger, and Wright
\cite{SWW1997}, who studied ordinary spherical means ($\alpha=0$) with
radii in a prescribed set $E\subset(0,\infty)$, obtaining sharp radial
estimates in terms of covering properties of $E$; for the
$L^p\to L^q$ problem on radial functions with fractal sets of radii, see
Beltran, Roos, and Seeger \cite{BRS2026}.

The endpoint behavior of the two segments is different.  More precisely,
in \cite[Corollary~1.6]{NRS2023} it was shown that if
$2\le p\le q_n$ and
$$
 \alpha=\frac1p-\frac{n-1}{2},
$$
then the radial operator is of restricted weak type $(p,p)$ but not of weak
type $(p,p)$, and hence is not of strong type.  By contrast, if $p>q_n$, the
boundary value
$$
 \alpha=-\frac{n-1}{p}
$$
is included in the strong-type radial range.  At the corner $p=q_n$, the
common boundary value belongs only to the restricted weak-type range.

The radial classification suggests the following conjecture for
general functions.

\begin{conjecture}[{\cite[Conjecture~1.7]{NRS2023}}]
\label{conj:NRS-general}
Let $n\ge2$, $\alpha\in\R$, and $1\le p<\infty$. Then $M^\alpha$ is
bounded on $L^p(\R^n)$ if and only if $p>1$ and
\begin{equation}
\label{eq:conjectured-general-region}
\begin{cases}
 \alpha>1-n+\dfrac np, &1<p<2,\\[5pt]
 \alpha>\dfrac1p-\dfrac{n-1}{2}, &2\le p\le q_n,\\[5pt]
 \alpha\ge-\dfrac{n-1}{p}, &p>q_n.
\end{cases}
\end{equation}
\end{conjecture}

\subsection{General necessary conditions}

Liu, Shen, Song, and Yan \cite{LSSY2024} subsequently completed the necessity picture for
general functions and complex orders.

\begin{theorem}[\cite{LSSY2024}]
\label{thm:LSSY}
Let $n\ge2$ and $2\le p<\infty$.  If $M^\alpha$ is bounded on $L^p(\R^n)$ for some
$\alpha\in\mathbb C$, then
\begin{equation}
\label{eq:general-necessity}
 \operatorname{Re}\alpha\ge
 \max\Big\{\frac1p-\frac{n-1}{2},-\frac{n-1}{p}\Big\}.
\end{equation}
\end{theorem}

For real $\alpha>(1-n)/2$, this necessary condition already follows from
Theorem~\ref{thm:NRS-radial} by restriction to radial functions.  In fact, on
the segment $2\le p\le q_n$, the radial theorem gives the slightly stronger
strict inequality.  Theorem~\ref{thm:LSSY} removes the auxiliary restriction
$\alpha>(1-n)/2$, applies to complex orders through
$\operatorname{Re}\alpha$, and in particular rules out boundedness whenever
$\operatorname{Re}\alpha\le(1-n)/2$.
Its proof gives a direct nonradial analysis of the two half-waves in the
asymptotic expansion of the multiplier and constructs examples adapted to
their possible interference.  An alternative proof was given by
F. Zhang \cite{Zhang2024}.

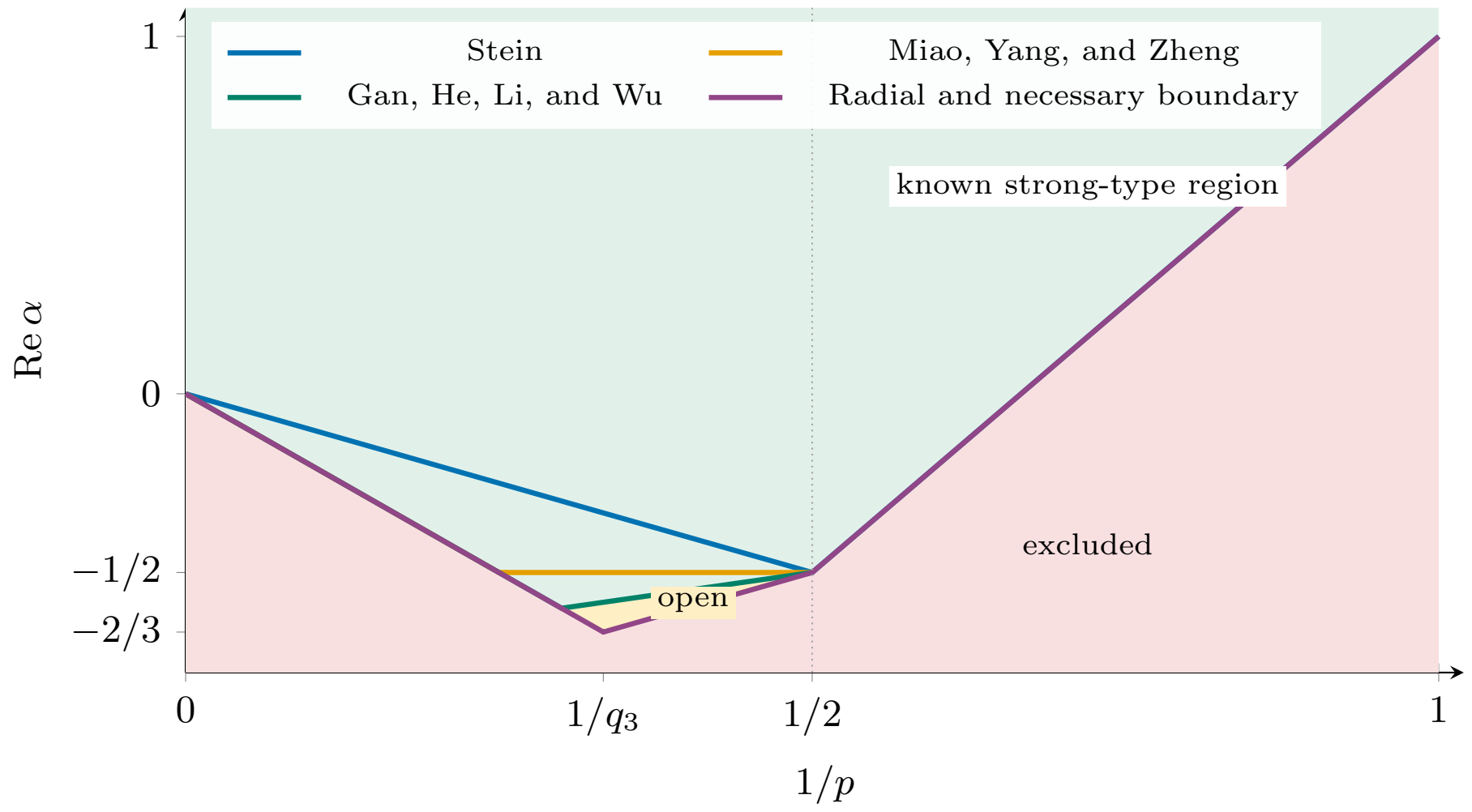
\begin{figure}[htbp]
\centering
\begin{tikzpicture}
\begin{axis}[
  width=0.72\textwidth,
  height=0.42\textwidth,
  xmin=0,
  xmax=1.02,
  ymin=-0.78,
  ymax=1.08,
  axis lines=left,
  xlabel={$1/p$},
  ylabel={$\operatorname{Re}\alpha$},
  xtick={0,0.333333,0.5,1},
  xticklabels={$0$,$1/q_3$,$1/2$,$1$},
  ytick={-0.666667,-0.5,0,1},
  yticklabels={$-2/3$,$-1/2$,$0$,$1$},
  tick label style={font=\small},
  label style={font=\small},
  legend style={
    at={(0.02,0.98)},
    anchor=north west,
    draw=none,
    fill=white,
    fill opacity=0.92,
    text opacity=1,
    font=\scriptsize,
    legend columns=2,
    column sep=8pt
  },
  clip=false
]
\path[draw=none,fill=ExcludedFill]
  (axis cs:0,-0.78) -- (axis cs:0,0)
  -- (axis cs:0.333333,-0.666667) -- (axis cs:0.5,-0.5)
  -- (axis cs:1,1) -- (axis cs:1,-0.78) -- cycle;
\path[draw=none,fill=GapFill]
  (axis cs:0.3,-0.6) -- (axis cs:0.5,-0.5)
  -- (axis cs:0.333333,-0.666667) -- cycle;
\path[draw=none,fill=KnownFill]
  (axis cs:0,1.08) -- (axis cs:1,1.08) -- (axis cs:1,1)
  -- (axis cs:0.5,-0.5) -- (axis cs:0.3,-0.6)
  -- (axis cs:0,0) -- cycle;

\addplot[SteinBlue,very thick]
  coordinates {(0,0) (0.5,-0.5) (1,1)};
\addlegendentry{Stein}

\addplot[MYZOrange,very thick]
  coordinates {(0,0) (0.25,-0.5) (0.5,-0.5)};
\addlegendentry{Miao, Yang, and Zheng}

\addplot[GHLWGreen,very thick]
  coordinates {(0,0) (0.3,-0.6) (0.5,-0.5)};
\addlegendentry{Gan, He, Li, and Wu}

\addplot[BoundaryPurple,very thick]
  coordinates {(0,0) (0.333333,-0.666667) (0.5,-0.5) (1,1)};
\addlegendentry{Radial and necessary boundary}

\addplot[forget plot,black!35,thin,densely dotted]
  coordinates {(0.5,-0.78) (0.5,1.08)};

\node[font=\scriptsize,fill=white,inner sep=1.5pt]
  at (axis cs:0.72,0.58) {known strong-type region};
\node[font=\scriptsize,fill=GapFill,inner sep=1.5pt]
  at (axis cs:0.405,-0.585) {open};
\node[font=\scriptsize,fill=ExcludedFill,inner sep=1.5pt]
  at (axis cs:0.72,-0.42) {excluded};
\end{axis}
\end{tikzpicture}
\caption{The boundedness diagram for $n=3$.  The colored segments are the
successive lower boundaries furnished by Stein, by Miao, Yang, and Zheng, see Section \ref{sec:BD-toMYZ}, and
as a consequence of the improved local smoothing by  Gan, He, Li, and Wu, see Section \ref{sec:GHLW}.  The purple curve is the radial and general necessary
boundary.  Green shading indicates the currently known open strong-type
region, yellow the remaining interior gap, and red the region excluded by
necessity.  The segments represent strict inequalities for positive results;
behavior at equality should be discussed separately.}
\label{fig:complete-phase-n3}
\end{figure}

The two quantities in \eqref{eq:general-necessity} agree at $p=q_n$.  It is
convenient to encode the conjectural open region by
\begin{equation}
\label{eq:necessary-epsilon}
 E_{\rm nec}(p,n)
 :=\min\Big\{\frac{n-1}{2}-\frac1p,\frac{n-1}{p}\Big\},
 \qquad
 \operatorname{Re}\alpha>-E_{\rm nec}(p,n).
\end{equation}
Figure~\ref{fig:complete-phase-n3} displays the whole boundedness diagram in the
representative case $n=3$.  The picture includes the already sharp range $1<p\le2$ so
that the improvement of the positive theory for $p>2$ can be seen in its
proper context, see Sections~\ref{sec:BD-toMYZ} and \ref{sec:GHLW}.

\section{From Bessel asymptotics to half-waves}
\label{sec:dictionary}

The connection between Stein's generalized spherical maximal operator and
the wave equation is already visible in the Fourier multiplier of
$\M_t^\alpha$.  We make it explicit before stating the transfer principle
in Section~\ref{sec:transfer}.

\subsection{The high-frequency Bessel decomposition}

Recall that 
$$
 \widehat{m_\alpha}(\xi)
 =c_{n,\alpha}|\xi|^{-\nu}J_\nu(|\xi|),
 \qquad
 \nu=\frac n2+\alpha-1,
$$
where $c_{n,\alpha}$ is given in \eqref{eq:multiplier}.
For fixed $\nu\in\mathbb C$, $r\ge1$, and an integer $N\ge1$, the Bessel
function has the expansion (see \cite[\S10.17(i), (ii)]{handbook})
\begin{equation}
\label{eq:Bessel-expansion}
 J_\nu(r)
 =r^{-1/2}\sum_{\varepsilon=\pm1}e^{i\varepsilon r}
   \Big(\sum_{k=0}^{N-1}c_{\varepsilon,k}(\nu)r^{-k}\Big)
   +R_{N,\nu}(r),
\end{equation}
where the remainder has the more precise form
\begin{equation}
\label{eq:Bessel-remainder}
 R_{N,\nu}(r)
=r^{-1/2}\sum_{\varepsilon=\pm1}
       e^{i\varepsilon r}\rho_{\varepsilon,N,\nu}(r),\qquad
 \big|\partial_r^\ell\rho_{\varepsilon,N,\nu}(r)\big|
 \lesssim_{N,\ell,\nu}r^{-N-\ell},\qquad \ell=0,1,2,\ldots.
\end{equation}
Thus the amplitudes $\rho_{\varepsilon,N,\nu}$ are symbols of order $-N$,
whereas $R_{N,\nu}$ satisfies
$|\partial_r^\ell R_{N,\nu}(r)|\lesssim_{N,\ell,\nu}r^{-N-1/2}$, since
derivatives can fall on $e^{i\varepsilon r}$.

For completeness, if
$$
 a_0(\nu)=1,
 \qquad
 a_k(\nu)=\frac{\prod_{m=1}^k
   \big(4\nu^2-(2m-1)^2\big)}{k!\,8^k},
$$
then the coefficients $c_{\varepsilon,k}(\nu)$ are given by
\begin{equation}
\label{eq:Bessel-coefficients}
 c_{\varepsilon,k}(\nu)
 =\frac{e^{-i\varepsilon(\pi\nu/2+\pi/4)}}{\sqrt{2\pi}}
   (i\varepsilon)^k a_k(\nu),
 \qquad \varepsilon=\pm1.
\end{equation}
Consequently, if
\begin{equation}
\label{eq:kappalpha}
 \kappa_\alpha:=\frac{n-1}{2}+\alpha,
\end{equation}
then, for $|\xi|\ge1$,
\begin{equation}
\label{eq:multiplier-half-wave}
 \widehat{m_\alpha}(\xi)
 =\sum_{\varepsilon=\pm1}e^{i\varepsilon|\xi|}
   |\xi|^{-\kappa_\alpha}b_{\varepsilon,\alpha}(|\xi|)
   +r_{N,\alpha}(|\xi|).
\end{equation}
Here, 
$$
 b_{\varepsilon,\alpha}(r)
 :=c_{n,\alpha}\sum_{k=0}^{N-1}c_{\varepsilon,k}(\nu)r^{-k},
 \qquad
 r_{N,\alpha}(r):=c_{n,\alpha}r^{-\nu}R_{N,\nu}(r).
$$
The functions $b_{\varepsilon,\alpha}$ are symbols of order zero on
$r\ge1$.  By \eqref{eq:Bessel-remainder}, $r_{N,\alpha}$ is a sum of the
two oscillations $e^{\pm ir}$ with amplitudes of order
$-\operatorname{Re}\kappa_\alpha-N$.  The symbol seminorms are locally uniform in $\alpha$ away from the
exceptional values of the normalized family.

Let $P_j$ be a Littlewood--Paley projection to
$|\xi|\simeq2^j$, and let $|D|=\sqrt{-\Delta}$ denote the Fourier
multiplier with symbol $|\xi|$.  For $j\ge1$ and $1\le t\le2$, formula
\eqref{eq:multiplier-half-wave} gives
\begin{equation}
\label{eq:dyadic-half-wave}
 \M_t^\alpha P_j f
 =\sum_{\varepsilon=\pm1}
 e^{i\varepsilon t\sqrt{-\Delta}}
 A_{\varepsilon,t}^\alpha
 |D|^{-\kappa_\alpha}P_jf
 +\mathcal R_{j,t}^\alpha f,
\end{equation}
where, on the support of $P_j$, the multiplier of
$A_{\varepsilon,t}^\alpha$ is
$t^{-\kappa_\alpha}b_{\varepsilon,\alpha}(t|\xi|)$, and
$$
 \widehat{\mathcal R_{j,t}^\alpha f}(\xi)
 =r_{N,\alpha}(t|\xi|)\widehat{P_jf}(\xi).
$$
The order-zero multiplier bounds are uniform for $t\in[1,2]$.
By taking sufficiently many terms in
\eqref{eq:Bessel-expansion}, the remainder $\mathcal R_{j,t}^\alpha$ has
any prescribed decay in $j$, also after any fixed finite number of
derivatives in $t$.  The low-frequency contribution comes from a smooth multiplier and is
handled by standard maximal-function estimates.  Thus the essential part of
$\M_t^\alpha$ is a sum of the two half-wave propagators, preceded by an
operator supplying
$$
 \operatorname{Re}\kappa_\alpha
 =\frac{n-1}{2}+\operatorname{Re}\alpha
$$
derivatives when this quantity is positive (and losing derivatives when it
is negative).
For positive estimates the two signs can be handled separately.  For
necessity, however, one must account for possible cancellation between the
incoming and outgoing waves and this is the interference issue
resolved by the examples of Liu, Shen, Song, and Yan in
Theorem~\ref{thm:LSSY}.

Miao, Yang, and Zheng \cite{MYZ2017} used Schl\"afli's integral
representation to carry out the decomposition in
\eqref{eq:dyadic-half-wave}.  For $r>0$ and $\nu\in\mathbb C$, it reads
(see also \cite[\S10.9.6]{handbook})
\begin{equation}
\label{eq:Schlafli}
 J_\nu(r)
 =\frac1{2\pi}\int_{-\pi}^{\pi}
   e^{i(r\sin\theta-\nu\theta)}\,d\theta
 -\frac{\sin(\pi\nu)}{\pi}\int_0^\infty
   e^{-r\sinh s-\nu s}\,ds.
\end{equation}
Treating $e^{-i\nu\theta}$ as part of the amplitude, the leading phase
$r\sin\theta$ in the first integral has stationary points at
$\theta=\pm\pi/2$, which produce the two oscillations $e^{\pm ir}$ in
\eqref{eq:Bessel-expansion}.  The second integral is nonoscillatory,
but it is not by itself the remainder in \eqref{eq:multiplier-half-wave}:
its expansion in powers of $r^{-1}$ is cancelled by the nonstationary
endpoint contributions of the first integral.

\subsection{Cone geometry and local smoothing}
\label{sub:conels}

For $g\in\mathcal S(\R^n)$, let
$$
 u_\varepsilon(x,t)=e^{i\varepsilon t\sqrt{-\Delta}}g(x),
 \qquad \varepsilon=\pm1.
$$
Without a time cutoff, the space-time Fourier transform of
$u_\varepsilon$ is supported on the sheet
$$
 \mathcal C_\varepsilon
 =\{(\xi,\tau)\in\R^n\times\R:\tau=\varepsilon|\xi|\}
$$
of the light cone.  Multiplication by a smooth time cutoff replaces this
support condition by rapid decay away from the cone in the temporal
frequency variable.  This is the geometric reason that estimates for
$M^\alpha$ are connected with cone restriction, decoupling, and square
function estimates.  A decomposition of the annulus into small angular caps
induces a decomposition of $u_\varepsilon$ into wave packets.  In physical
space, these packets are concentrated on tubes whose directions are normal
to the corresponding caps of the cone.  Decoupling inequalities measure the
extent to which packets with different directions can be treated almost
orthogonally, while reverse square function estimates reassemble the angular
pieces without losing the geometry of their interactions.

For $s\in\R$ and $1<p<\infty$, we use the inhomogeneous Bessel-potential
Sobolev space
\begin{equation}
\label{eq:Sobolev-definition}
 W^{s,p}(\R^n)
 :=\{f\in\mathcal S'(\R^n):\langle D\rangle^s f\in L^p(\R^n)\}, \qquad
 \|f\|_{W^{s,p}}:=\|\langle D\rangle^s f\|_{L^p},
 \qquad \langle D\rangle^s:=(I-\Delta)^{s/2}.
\end{equation}
Thus $\widehat{\langle D\rangle^sf}(\xi)
=(1+|\xi|^2)^{s/2}\widehat f(\xi)$.  This specifies the convention also
when $s$ is nonintegral or negative. For later use, write
\begin{equation}
\label{eq:gammapn}
 \gamma(p,n):=(n-1)\Big(\frac12-\frac1p\Big),
 \qquad 2\le p<\infty.
\end{equation}
At a fixed time, the classical half-wave estimates of Miyachi and Peral
\cite{Miyachi1980,Peral1980} give
\begin{equation}
\label{eq:half-wave-fixed-time}
 \big\|e^{it|D|}f\big\|_{L^p(\R^n)}
 \lesssim \|f\|_{W^{\gamma(p,n),p}(\R^n)},
 \qquad 1\le t\le2,
\end{equation}
uniformly in $t$.
The general FIO version, see Subsection \ref{sub:FIO}, is due to Seeger, Sogge, and Stein
\cite{SSS1991}.
Local smoothing asks how much of this derivative loss can be recovered when
the solution is also averaged over a compact time interval.  Sogge's
Euclidean local smoothing conjecture \cite{Sogge1991} asserts that, for
$p\ge q_n$, with $q_n$ given in \eqref{eq:qn}, and every
\begin{equation}
\label{eq:slsc}
 s>\gamma(p,n)-\frac1p
 =\frac{n-1}{2}-\frac np,
\end{equation}
one has
\begin{equation}
\label{eq:Euclidean-LS-conjecture}
 \big\|e^{it|D|}f\big\|_{L^p(\R^n\times[1,2])}
 \lesssim_s \|f\|_{W^{s,p}(\R^n)}.
\end{equation}
For a systematic account and the formulation used here, see
\cite[\S2.3.1, Conjecture~20]{BHS2021}.
By Plancherel and the unitarity of $e^{it|D|}$,
\begin{equation}
\label{eq:half-wave-L2}
 \|e^{it|D|}f\|_{L^2(\R^n\times[1,2])}=\|f\|_{L^2(\R^n)}.
\end{equation}
The strict inequality for $s$ in \eqref{eq:slsc} allows an
arbitrarily small loss relative to the conjectured regularity
$s=\gamma(p,n)-1/p$.  Below $q_n$, interpolation with the $L^2$ identity predicts the
estimate for every $s>0$.  Thus the expected regularity over the whole range
$2<p<\infty$ is
$$
 s>s_{\rm Euc}(p,n):=
 \max\Big\{0,\frac{n-1}{2}-\frac np\Big\}.
$$

The Sobolev regularity exponent $s$ in \eqref{eq:Euclidean-LS-conjecture} is the input used to obtain bounds for Stein's maximal operator.  Passing from an $L^p$ norm in
$(x,t)$ to a supremum in $t$ costs $1/p$ of a time derivative, while the
factor $|D|^{-\kappa_\alpha}$ in \eqref{eq:dyadic-half-wave} supplies
$(n-1)/2+\operatorname{Re}\alpha$ derivatives.  The precise calculation is
given in Section~\ref{sec:transfer}.

\subsection{The Fourier integral operator (FIO) formulation}
\label{sub:FIO}

The half-wave propagator is the simplest member of a broader class with
variable coefficients.  A local wave parametrix on a Riemannian manifold
has the form
\begin{equation}
\label{eq:FIO-model}
 \mathcal Ff(x,t)
 =\frac1{(2\pi)^n}\int_{\R^n}e^{i\phi(x,t;\xi)}
   a(x,t;\xi)\widehat f(\xi)\,d\xi,
\end{equation}
where $\phi$ is real-valued, smooth away from $\xi=0$, and positively
homogeneous of degree one in $\xi$.  An amplitude $a$ is a symbol of order
$\mu\in\R$, written $a\in S^\mu_{1,0}$, if
\begin{equation}
\label{eq:FIO-symbol-order}
 \big|\partial_{x,t}^{\beta}\partial_\xi^{\eta}a(x,t;\xi)\big|
 \le C_{\beta,\eta}(1+|\xi|)^{\mu-|\eta|}
\end{equation}
for all multiindices $\beta,\eta$.  Local FIO estimates use amplitudes
compactly supported in $(x,t)$. Throughout this survey, the order of
such a family means the order $\mu$ of its amplitude in
\eqref{eq:FIO-symbol-order}.  Thus $a=1$ and the phase
$$
 \phi(x,t;\xi)=x\cdot\xi+t|\xi|
$$
give the half-wave propagator $e^{it|D|}$ as an operator of order
zero, and multiplying the amplitude by $|\xi|^{-\kappa_\alpha}$ at high
frequency gives order $-\operatorname{Re}\kappa_\alpha$.  In H\"ormander's
convention \cite[\S3.2]{Hormander1971}, the fixed-time operators have the
same order $\mu$, whereas the space--time operator belongs to the class
$I^{\mu-1/4}$; see also \cite{BHS2021}.

The relevant phases satisfy two geometric hypotheses.  Writing $z=(x,t)$,
the nondegeneracy condition is
$$
 \operatorname{rank}\partial_{z\xi}^2\phi(z;\xi)=n.
$$
In wave-equation charts one also has
$\det\partial_{x\xi}^2\phi\ne0$, which gives a canonical graph at each
fixed time.  The curvature condition concerns the conic hypersurface
$$
 \Sigma_z:=\{\partial_z\phi(z;\xi):\xi\ne0\}
 \subset T_z^*\R^{n+1},
$$
restricted to the frequency directions under consideration.  It must have
$n-1$ nonzero principal curvatures; the radial direction has zero curvature
by homogeneity.  In the positive-definite case, the nonzero curvatures have
the same sign.  These are the hypotheses of cinematic-curvature theory;
see Mockenhaupt, Seeger, and Sogge \cite{MSS1993}, the survey of Beltran,
Hickman, and Sogge \cite{BHS2021}, and Sogge's monograph
\cite[Ch.~8]{Sogge2017}.

After the Bessel decomposition, the principal terms in
$\M_t^\alpha$ are FIO families of amplitude order
$$
 -\frac{n-1}{2}-\operatorname{Re}\alpha
$$
with the Euclidean half-wave phase.  Corresponding local smoothing theorems
are available for variable-coefficient phases satisfying the stated
hypotheses, but their admissible range of $p$ need not coincide with the
Euclidean one; see Section~\ref{sec:GHLW}.  This places Stein's family
naturally within FIO theory and suggests a variable-coefficient analogue;
see Section~\ref{sec:open}.

\section{The transfer principle}
\label{sec:transfer}

We now explain the mechanism that converts a local smoothing estimate into a
bound for Stein's maximal operator.  The reduction goes back to \cite{MYZ2017}; the formulation below, which keeps both competing
exponents visible, is essentially \cite[Proposition~4.1]{LSSY2024}.  Let
$2<p<\infty$, and assume that, for
some $s\in\R$,
\begin{equation}
\label{eq:LSs}
 \|e^{it|D|}f\|_{L^p(\R^n\times[1,2])}
 \lesssim \|f\|_{W^{s,p}(\R^n)}.
\end{equation}
The estimate for $e^{-it|D|}$ follows by complex conjugation.  In particular,
if $P_j$ is supported where $|\xi|\simeq2^j$, assumption
\eqref{eq:LSs} implies
\begin{equation}
\label{eq:localized-LSs}
 \|e^{\pm it|D|}P_jf\|_{L^p(\R^n\times[1,2])}
 \lesssim 2^{js}\|P_jf\|_{L^p(\R^n)},
 \qquad j\ge1.
\end{equation}
The estimate remains valid after inserting symbols
$a_j(t,\xi)$ with uniform order-zero bounds, supported in a slightly larger
annulus.  Indeed, for each $t$ the corresponding Fourier multiplier is
uniformly bounded on $L^p_x$ and commutes with $e^{\pm it|D|}$.

\begin{proposition}[Transfer principle]
\label{prop:transfer}
Suppose that \eqref{eq:LSs} holds.  Then
\begin{equation}
\label{eq:transfer-condition}
 \|M^\alpha f\|_{L^p(\R^n)}\lesssim\|f\|_{L^p(\R^n)}
 \quad\text{provided}\quad
 \operatorname{Re}\alpha>
 \max\Big\{-\frac{n-1}{p},
 s-\frac{n-1}{2}+\frac1p\Big\}.
\end{equation}
Equivalently, if
\begin{equation}
\label{eq:transfer-Es}
 E_s(p,n):=\min\Big\{\frac{n-1}{p},
 \frac{n-1}{2}-\frac1p-s\Big\},
\end{equation}
then $\operatorname{Re}\alpha>-E_s(p,n)$ is sufficient.
\end{proposition}

We will give the proof below; every later application amounts to substituting a value of $s$ into \eqref{eq:transfer-condition}.  The following elementary
one-dimensional estimate is the device that converts space--time control into
a supremum in $t$.

\begin{lemma}
\label{lem:time-maximal}
Let $1<p<\infty$ and let $F$ be smooth on
$\R^n\times[1,2]$.  Then
\begin{equation}
\label{eq:time-maximal}
 \big\|\sup_{1\le t\le2}|F(\cdot,t)|\big\|_{L^p_x}
 \lesssim_p
 \|F(\cdot,1)\|_{L^p_x}
 +\|F\|_{L^p_{x,t}}^{1-1/p}
  \|\partial_tF\|_{L^p_{x,t}}^{1/p}.
\end{equation}
\end{lemma}
\begin{proof}
The fundamental theorem of calculus gives, pointwise in $x$,
$$
 \sup_{1\le t\le2}|F(x,t)|^p
 \le |F(x,1)|^p+
   p\int_1^2|F(x,t)|^{p-1}|\partial_tF(x,t)|\,dt.
$$
Integration in $x$ and H\"older's inequality on the product space give
$$
 \int_{\R^n}\int_1^2
 |F(x,t)|^{p-1}|\partial_tF(x,t)|\,dt\,dx
 \le \|F\|_{L^p_{x,t}}^{p-1}
       \|\partial_tF\|_{L^p_{x,t}},
$$
which proves \eqref{eq:time-maximal} after taking the $p$th root.  One can
equivalently apply H\"older first in $t$ and then in $x$, or in the reverse
order.
\end{proof}

If a frequency-localized family satisfies
$\|F\|_{L^p_{x,t}}\lesssim A_j$ and
$\|\partial_tF\|_{L^p_{x,t}}\lesssim2^j A_j$, the mixed term in
\eqref{eq:time-maximal} is at most $2^{j/p}A_j$. Miao, Yang, and
Zheng \cite{MYZ2017} use instead the Sobolev embedding
$W^{1/p+\eta,p}([1,2])\hookrightarrow L^\infty([1,2])$, $\eta>0$, which
costs $2^{j(1/p+\eta)}$; since $\eta$ is absorbed by the strict inequality
in $\alpha$, both approaches give the same open region.

\begin{proof}[Proof of Proposition~\ref{prop:transfer}]
Choose $\chi_0\in C_c^\infty([0,\infty))$ equal to $1$ on $[0,1]$ and
to $0$ on $[2,\infty)$, and set $\chi(r)=\chi_0(r)-\chi_0(2r)$.  Then
$$
 \chi_0(r)+\sum_{j\ge1}\chi_j(r)=1,
 \qquad \chi_j(r):=\chi(2^{-j}r),
$$
with $\chi_j$ supported in $[2^{j-1},2^{j+1}]$ for $j\ge1$.  The
fact that $\chi_0$ is constant near zero ensures that $\chi_0(|\xi|)$ is smooth.
The decomposition is made
in the scale-invariant variable $t|\xi|$.  Thus define
\begin{equation}
\label{eq:relative-frequency-piece}
 \widehat{\M_{j,t}^\alpha f}(\xi)
 :=\chi_j(t|\xi|)\widehat m_\alpha(t\xi)\widehat f(\xi),
 \qquad j\ge0.
\end{equation}

For $j=0$, set
$$
 \eta_0^\alpha(\xi)
 :=\chi_0(|\xi|)\widehat m_\alpha(\xi).
$$
This is a smooth compactly supported multiplier.  Its inverse Fourier
transform $K_0^\alpha$ is a Schwartz function, so there is a radial,
decreasing function $\Phi_\alpha\in L^1(\R^n)$ such that
$|K_0^\alpha(x)|\le\Phi_\alpha(x)$.  Since the multiplier in
\eqref{eq:relative-frequency-piece} is
$\eta_0^\alpha(t\xi)$, its kernel is the usual dilation
$K_{0,t}^\alpha(x)=t^{-n}K_0^\alpha(x/t)$.  The standard domination of
convolution with dilates of an integrable radial decreasing function by the
Hardy--Littlewood maximal operator therefore gives
\begin{equation}
\label{eq:transfer-low-frequency}
 \sup_{t>0}|\M_{0,t}^\alpha f(x)|
 \lesssim_\alpha M_{\rm HL}f(x).
\end{equation}
It remains to treat $j\ge1$.

Recall $\kappa_\alpha$ from \eqref{eq:kappalpha}.  Fix the truncation order
$N$ in \eqref{eq:Bessel-expansion}; it will be chosen sufficiently large
below.  Apply \eqref{eq:multiplier-half-wave} with
$r=t|\xi|$ on the support of $\chi_j(t|\xi|)$.  This gives the identity
\begin{equation}
\label{eq:transfer-half-wave-piece}
 \M_{j,t}^\alpha f
 =2^{-j\kappa_\alpha}
   \sum_{\varepsilon=\pm1}T_{\varepsilon,j}^\alpha f(\cdot,t)
   +\mathcal E_{j,t}^\alpha f.
\end{equation}
Here
$$
 T_{\varepsilon,j}^\alpha f(x,t)
 =\frac1{(2\pi)^n}\int_{\R^n}
   e^{i(x\cdot\xi+\varepsilon t|\xi|)}
   a_{\varepsilon,j}^\alpha(t,\xi)\widehat f(\xi)\,d\xi,
$$
with the amplitude
$$
 a_{\varepsilon,j}^\alpha(t,\xi)
 :=\chi_j(t|\xi|)(2^{-j}t|\xi|)^{-\kappa_\alpha}
       b_{\varepsilon,\alpha}(t|\xi|).
$$
It is supported where $|\xi|\simeq2^j$ when $1\le t\le2$, and satisfies
$$
 \big|\partial_t^k\partial_\xi^\beta
     a_{\varepsilon,j}^\alpha(t,\xi)\big|
 \lesssim_{k,\beta,N,\alpha}2^{-j|\beta|}.
$$
The remainder in \eqref{eq:transfer-half-wave-piece} is defined by
\begin{equation}
\label{eq:transfer-remainder-definition}
 \widehat{\mathcal E_{j,t}^\alpha f}(\xi)
 :=\chi_j(t|\xi|)r_{N,\alpha}(t|\xi|)\widehat f(\xi).
\end{equation}
It uses the same scalar remainder $r_{N,\alpha}$ as
$\mathcal R_{j,t}^\alpha$ in \eqref{eq:dyadic-half-wave}; the distinction
is the relative-frequency cutoff $\chi_j(t|\xi|)$ in place of $P_j$.
The dependence on $N$ is suppressed in both operator notations.

By \eqref{eq:Bessel-remainder}, this remainder is a sum of two half-waves
whose amplitudes have order $-\operatorname{Re}\kappa_\alpha-N$.
The fixed-time estimate \eqref{eq:half-wave-fixed-time}, followed by
Lemma~\ref{lem:time-maximal}, therefore gives
$$
 \big\|\sup_{1\le t\le2}|\mathcal E_{j,t}^\alpha f|\big\|_{L^p}
 \lesssim_{N,\alpha}
 2^{-j(N+\operatorname{Re}\alpha+(n-2)/p)}
 \|\widetilde P_jf\|_{L^p},
$$
where $\widetilde P_j$ is a slightly enlarged Littlewood--Paley projection.
Indeed, before taking the time supremum the exponent is
$-\operatorname{Re}\kappa_\alpha-N+\gamma(p,n)$, with $\gamma(p,n)$ defined in \eqref{eq:gammapn}; one time derivative
costs at most $2^j$, and the lemma costs $2^{j/p}$. The truncation order
$N$ can be arbitrarily large. Given any $N_0>0$, choose it so that
$N+\operatorname{Re}\alpha+(n-2)/p\ge N_0$. Then
\begin{equation}
\label{eq:transfer-remainder}
 \big\|\sup_{1\le t\le2}|\mathcal E_{j,t}^\alpha f|\big\|_{L^p}
 \lesssim_{N_0,\alpha}2^{-jN_0}\|\widetilde P_jf\|_{L^p}.
\end{equation}

We next estimate the oscillatory terms. Applying
\eqref{eq:half-wave-fixed-time} to each half-wave in
\eqref{eq:transfer-half-wave-piece}, with its uniformly bounded
order-zero amplitude, gives
\begin{equation}
\label{eq:transfer-fixed-time}
 \|T_{\varepsilon,j}^\alpha f(\cdot,1)\|_{L^p}
 \lesssim
 2^{j\gamma(p,n)}\|\widetilde P_jf\|_{L^p};
\end{equation}
For the space-time bounds, we use the frequency-localized local smoothing
estimate \eqref{eq:localized-LSs}. Its stability under uniformly bounded
order-zero amplitudes yields
\begin{align}
\label{eq:transfer-LS-piece}
 \|T_{\varepsilon,j}^\alpha f\|_{L^p_{x,t}}
 &\lesssim 2^{js}\|\widetilde P_jf\|_{L^p},\\
\label{eq:transfer-LS-derivative}
 \|\partial_tT_{\varepsilon,j}^\alpha f\|_{L^p_{x,t}}
 &\lesssim 2^{j(s+1)}\|\widetilde P_jf\|_{L^p}.
\end{align}
The extra factor $2^j$ in \eqref{eq:transfer-LS-derivative} comes from
differentiating the phase.  
Lemma~\ref{lem:time-maximal} therefore gives
\begin{equation}
\label{eq:transfer-main-maximal}
 \big\|\sup_{1\le t\le2}
   |T_{\varepsilon,j}^\alpha f(\cdot,t)|\big\|_{L^p}
 \lesssim
 \big(2^{j\gamma(p,n)}+2^{j(s+1/p)}\big)
 \|\widetilde P_jf\|_{L^p}.
\end{equation}

Since $|2^{-j\kappa_\alpha}|
=2^{-j((n-1)/2+\operatorname{Re}\alpha)}$, multiplication by the
first term on the right of \eqref{eq:transfer-main-maximal}, which
comes from the value at $t=1$ in Lemma~\ref{lem:time-maximal}, gives
$$
 2^{-j((n-1)/2+\operatorname{Re}\alpha)}
 2^{j\gamma(p,n)}
 =
 2^{-j(\operatorname{Re}\alpha+(n-1)/p)}.
$$
The local smoothing term gives
$$
 2^{-j((n-1)/2+\operatorname{Re}\alpha)}
 2^{j(s+1/p)}
 =
 2^{-j(\operatorname{Re}\alpha+(n-1)/2-s-1/p)}.
$$
Thus the hypothesis in \eqref{eq:transfer-condition} says that both
decay exponents are positive.  Taking
$$
 0<\delta<
 \min\Big\{
 \operatorname{Re}\alpha+\frac{n-1}{p},
 \operatorname{Re}\alpha+\frac{n-1}{2}-s-\frac1p
 \Big\}
$$
and choosing $N_0$ large enough in \eqref{eq:transfer-remainder}, we obtain
the local-in-radius estimate
\begin{equation}
\label{eq:transfer-local-radius}
 \big\|\sup_{1\le t\le2}|\M_{j,t}^\alpha f|\big\|_{L^p}
 \lesssim 2^{-j\delta}\|\widetilde P_jf\|_{L^p},
 \qquad j\ge1.
\end{equation}

We now pass from $1\le t\le2$ to all radii. Let $I_k=[2^k,2^{k+1}]$. If
$t\in I_k$ and $t|\xi|\simeq2^j$, then $|\xi|\simeq2^{j-k}$.
Let $\widetilde P_\ell$, $\ell\in\mathbb Z$, now denote homogeneous
Littlewood--Paley projections, with the same enlargement  as
above.  Write $t=2^k\tau$ and $x=2^ky$, and set
$f_k(y)=f(2^ky)$.  The identity
$$
 \M_{2^k\tau}^\alpha f(2^ky)=\M_\tau^\alpha f_k(y)
$$
reduces $t\in I_k$ to $\tau\in[1,2]$.  Moreover, a frequency
$2^{j-k}$ for $f$ becomes a frequency $2^j$ for $f_k$, and the scaling
factors in the two $L^p$ norms cancel.  Hence
\eqref{eq:transfer-local-radius} implies
\begin{equation}
\label{eq:transfer-scaled-radius}
 \big\|\sup_{t\in I_k}|\M_{j,t}^\alpha f|\big\|_{L^p}
 \lesssim 2^{-j\delta}\|\widetilde P_{j-k}f\|_{L^p},
\end{equation}
uniformly in $k\in\mathbb Z$.  Consequently,
\begin{align*}
 \big\|\sup_{t>0}|\M_{j,t}^\alpha f|\big\|_{L^p}
 &\le
 \Big(\sum_{k\in\mathbb Z}
  \big\|\sup_{t\in I_k}|\M_{j,t}^\alpha f|\big\|_{L^p}^p
 \Big)^{1/p}
 \lesssim 2^{-j\delta}
 \Big\|\Big(\sum_{k\in\mathbb Z}
 |\widetilde P_{j-k}f|^p\Big)^{1/p}\Big\|_{L^p}\\
 &\le 2^{-j\delta}
 \Big\|\Big(\sum_{k\in\mathbb Z}
 |\widetilde P_{j-k}f|^2\Big)^{1/2}\Big\|_{L^p}
 \lesssim 2^{-j\delta}\|f\|_{L^p}.  
\end{align*}
The penultimate inequality uses $p\ge2$, as assumed here.  Summing the
geometric series over $j\ge1$ and combining it with
\eqref{eq:transfer-low-frequency} proves the proposition.
\end{proof}

The two terms in \eqref{eq:transfer-condition} come, respectively, from
the fixed-time estimate \eqref{eq:transfer-fixed-time} and local smoothing
with the $1/p$ cost of the time supremum.  At the conjectured local
smoothing index $s_{\rm Euc}(p,n)$, they coincide for
$p\ge q_n=2n/(n-1)$.

\begin{remark}[Strict inequalities]
\label{rem:strict-inequalities}
Local smoothing theorems are usually stated with an arbitrarily small loss:
\eqref{eq:LSs} holds for every $s>s_0$.  Given
$$
 \operatorname{Re}\alpha>
 \max\Big\{-\frac{n-1}{p},
 s_0-\frac{n-1}{2}+\frac1p\Big\},
$$
one first chooses $s>s_0$ sufficiently close to $s_0$ and then chooses
$\delta>0$ above.  The strict inequality both absorbs the local smoothing
loss and ensures a summable geometric series in $j$.  Even a local smoothing
estimate at $s=s_0$ would not settle equality in
$\operatorname{Re}\alpha$.  Boundary estimates require a separate argument.
\end{remark}

\section{From Bourgain and Demeter to the theorem of Miao, Yang, and Zheng}
\label{sec:BD-toMYZ}

The first application of the transfer principle came from the cone
decoupling theorem of Bourgain and Demeter \cite{BD2015}.  At frequency
$|\xi|\simeq2^j$, their theorem implies the sharp conjectural local smoothing
estimate in the range
$$
 p\ge p_{\rm BD}(n):=\frac{2(n+1)}{n-1}.
$$
In the notation of \eqref{eq:LSs}, this means that one may take every
\begin{equation}
\label{eq:BD-large-p}
 s>\frac{n-1}{2}-\frac np,
 \qquad p\ge p_{\rm BD}(n).
\end{equation}
For $2<p<p_{\rm BD}(n)$, interpolating the estimate at
$p=p_{\rm BD}(n)$ with the $L^2$ identity \eqref{eq:half-wave-L2} gives
\begin{equation}
\label{eq:BD-interpolated}
 s>\frac{n-1}{2}\Big(\frac12-\frac1p\Big).
\end{equation}
Thus the input of Bourgain and Demeter, together with interpolation, can be
recorded for all $2<p<\infty$ as
\begin{equation}
\label{eq:BDs}
 s>s_{\rm BD}(p,n),
 \qquad
 s_{\rm BD}(p,n):=
 \max\Big\{
 \frac{n-1}{2}\Big(\frac12-\frac1p\Big),
 \frac{n-1}{2}-\frac np
 \Big\}.
\end{equation}
The two expressions in the maximum in \eqref{eq:BDs} agree at
$p=p_{\rm BD}(n)$.  The first one is therefore the regularity obtained by
interpolation, whereas the second is the sharp
local smoothing consequence of cone decoupling.  This distinction will be
important when the input of Bourgain and Demeter is replaced by later local
smoothing theorems with a smaller threshold.

We now insert \eqref{eq:BDs} into Proposition~\ref{prop:transfer}. By
\eqref{eq:transfer-Es}, the admissible negative order supplied by local
smoothing with exponent $s$ is $E_s(p,n)$.
In the interpolation range $2<p<p_{\rm BD}(n)$, a direct calculation gives
\begin{align}
\label{eq:MYZ-interpolation-calculation}
 \frac{n-1}{2}-\frac1p
 -\frac{n-1}{2}\Big(\frac12-\frac1p\Big)
 &=\frac{n-1}{4}+\frac{n-3}{2p}.
\end{align}
On the decoupling branch $p\ge p_{\rm BD}(n)$, one instead obtains
\begin{align}
\label{eq:MYZ-decoupling-calculation}
 \frac{n-1}{2}-\frac1p
 -\Big(\frac{n-1}{2}-\frac np\Big)
 &=\frac{n-1}{p}.
\end{align}
Notice also that
\begin{equation}
\label{eq:MYZ-branch-comparison}
 \Big(\frac{n-1}{4}+\frac{n-3}{2p}\Big)
 -\frac{n-1}{p}
 =\frac{n-1}{4p}\big(p-p_{\rm BD}(n)\big).
\end{equation}
Consequently the first quantity is the smaller one below
$p_{\rm BD}(n)$, the second is the smaller one above it, and they agree at
the threshold.  We recover the following theorem.

\begin{theorem}[\cite{MYZ2017}]
\label{thm:MYZ}
Let $n\ge2$ and $2<p<\infty$.  Then $M^\alpha$ is bounded on
$L^p(\R^n)$ whenever
\begin{equation}
\label{eq:MYZepsilon}
 \operatorname{Re}\alpha>-E_{\rm MYZ}(p,n),
 \qquad
 E_{\rm MYZ}(p,n):=
 \min\Big\{
 \frac{n-1}{4}+\frac{n-3}{2p},
 \frac{n-1}{p}
 \Big\}.
\end{equation}
Equivalently,
\begin{equation}
\label{eq:MYZ-piecewise}
 \operatorname{Re}\alpha>
 \begin{cases}
 -\dfrac{n-1}{4}-\dfrac{n-3}{2p},
   &2<p<p_{\rm BD}(n),\\[6pt]
 -\dfrac{n-1}{p},
   &p\ge p_{\rm BD}(n).
 \end{cases}
\end{equation}
\end{theorem}

The first entry in the minimum in \eqref{eq:MYZepsilon} is the local
smoothing contribution on the interpolated branch, while the second is the
fixed-time FIO contribution.  At and above
$p_{\rm BD}(n)$ the two ranges balance in such a way that the transfer
argument reaches the necessary line
$\operatorname{Re}\alpha=-(n-1)/p$.  Thus Theorem \ref{thm:MYZ} is sharp, apart from the boundary, throughout
$$
 p\ge p_{\rm BD}(n).
$$

It is also instructive to compare this with Stein's classical positive
region.  For $p\ge2$, Stein's condition is
$$
 \operatorname{Re}\alpha>\frac{2-n}{p}.
$$
In the decoupling range, Miao, Yang, and Zheng shift this boundary downward
by $1/p$, obtaining $-(n-1)/p$.  Below
$p_{\rm BD}(n)$, the improvement is
\begin{align*}
 \frac{2-n}{p}
 -\Big(-\frac{n-1}{4}-\frac{n-3}{2p}\Big)
 &=\frac{n-1}{2}\Big(\frac12-\frac1p\Big),
\end{align*}
which is precisely the gain produced by interpolating the estimate of
Bourgain and Demeter with $L^2$.

\subsection{Comparison with the necessary region}

By \eqref{eq:MYZepsilon} and \eqref{eq:necessary-epsilon}, the two
quantities $E_{\rm MYZ}(p,n)$ and $E_{\rm nec}(p,n)$ coincide, and equal
$(n-1)/p$, for $p\ge p_{\rm BD}(n)$, whereas for $2<p<p_{\rm BD}(n)$ one has
$E_{\rm MYZ}(p,n)=\frac{n-1}{4}+\frac{n-3}{2p}<E_{\rm nec}(p,n)$.  The gap
has two different meanings.  For $2<p\le q_n$ the target is
$\operatorname{Re}\alpha>1/p-(n-1)/2$, and the difference
$$
 E_{\rm nec}(p,n)-E_{\rm MYZ}(p,n)
 =\frac{n-1}{2}\Big(\frac12-\frac1p\Big)
$$
is the price of interpolating from $p_{\rm BD}(n)$ instead of from $q_n$:
sharp local smoothing at $p=q_n$ would give every $s>0$ for $2<p<q_n$ and
close this part of the gap.  For $q_n<p<p_{\rm BD}(n)$ the target is
$\operatorname{Re}\alpha>-(n-1)/p$ and the difference is
$$
 E_{\rm nec}(p,n)-E_{\rm MYZ}(p,n)
 =\frac{n-1}{4p}\big(p_{\rm BD}(n)-p\big),
$$
a threshold gap that closes at $p_{\rm BD}(n)$.  Both differences are
largest at $p=q_n$, where they equal $(n-1)/(4n)$.  This explains why a
lower local smoothing threshold first shrinks the second interval and, only
after interpolation with the $L^2$ estimate, improves the first one.

\section{Dimension two: the complete open region}
\label{sec:n2}

Dimension two is the only dimension in which the local smoothing input
currently determines the entire open strong-type region.  In this case the
two distinguished exponents introduced above are
$$
 q_2=4
 \qquad\text{and}\qquad
 p_{\rm BD}(2)=6.
$$
Thus the input of Bourgain and Demeter already gives the conjectured maximal
estimate for $p\ge6$, but leaves two gaps: the interval $4\le p<6$, where
the conjectural boundary is $\operatorname{Re}\alpha=-1/p$, and the interval
$2<p<4$, where it is
$\operatorname{Re}\alpha=1/p-1/2$.  The theorem of Guth, H. Wang, and R. Zhang
closes the first gap directly and the second one by interpolation.

\subsection{The sharp planar local smoothing input}

The reverse square function estimate for the cone in $\R^3$ proved in \cite{GWZ2020}, combined with the standard implication from
cone square functions to local smoothing \cite{MSS1992}, gives the following
estimate for the half-wave propagator.

\begin{theorem}[Planar local smoothing, \cite{GWZ2020}]
\label{thm:GWZ-local-smoothing}
Let $4\le p<\infty$.  For every
$$
 s>\frac12-\frac2p
$$
one has
$$
 \big\|e^{it|D|}f\big\|_{L^p(\R^2\times[1,2])}
 \lesssim_s
 \|f\|_{W^{s,p}(\R^2)}.
$$
\end{theorem}

The range $2<p<4$ follows by interpolating the $p=4$ estimate with
the $L^2$ identity \eqref{eq:half-wave-L2}.
Indeed, if $2<p<4$, choose $0<\theta<1$ so that
$$
 \frac1p=\frac{1-\theta}{2}+\frac{\theta}{4}.
$$
Interpolating the $L^2$ estimate with the $L^4$ local smoothing estimate
having regularity $s_1>0$ gives regularity $\theta s_1$.  Given any $s>0$,
we may take $s_1=s/\theta$, and hence obtain the estimate with regularity
$s$.  It is therefore convenient to summarize the planar input as
$$
 s>s_{\rm pl}(p):=
 \max\Big\{0,\frac12-\frac2p\Big\}
 =
 \begin{cases}
  0, & 2<p\le4,\\
  \dfrac12-\dfrac2p, & 4\le p<\infty.
 \end{cases}
$$
Here and below the value at $p=4$ is  the corner
at which the two branches meet.

\subsection{Insertion into the transfer principle}

For $n=2$, Proposition~\ref{prop:transfer} reads
\begin{equation}
\label{eq:transfer-n2}
 \operatorname{Re}\alpha>
 \max\Big\{-\frac1p,s-\frac12+\frac1p\Big\}.
\end{equation}
Suppose first that $p\ge4$.  Taking $s$ arbitrarily close to
$1/2-2/p$ from above, the second term on the right-hand side of
\eqref{eq:transfer-n2} becomes
$$
 \Big(\frac12-\frac2p\Big)-\frac12+\frac1p
 =-\frac1p.
$$
Both terms in the maximum therefore agree, and the sufficient condition is
$\operatorname{Re}\alpha>-1/p$.

If $2<p<4$, we may instead take $s>0$ arbitrarily small.  The limiting
value of the second term is $1/p-1/2$, and
$$
 \Big(\frac1p-\frac12\Big)-\Big(-\frac1p\Big)
 =\frac2p-\frac12>0.
$$
Consequently the second term controls the maximum and the sufficient
condition is $\operatorname{Re}\alpha>1/p-1/2$.  At $p=4$ the two
expressions coincide and equal $-1/4$.  We have thus obtained the complete
open region.

\begin{theorem}[The planar open region]
\label{thm:n2-open-region}
Let $2<p<\infty$ and $\alpha\in\mathbb C$.
\begin{enumerate}
\item If
\begin{equation}
\label{eq:n2sharp}
 \operatorname{Re}\alpha>
 \max\Big\{\frac1p-\frac12,-\frac1p\Big\},
\end{equation}
then $M^\alpha$ is bounded on $L^p(\R^2)$.
\item If $M^\alpha$ is bounded on $L^p(\R^2)$, then
\begin{equation}
\label{eq:n2-necessary}
 \operatorname{Re}\alpha\ge
 \max\Big\{\frac1p-\frac12,-\frac1p\Big\}.
\end{equation}
\end{enumerate}
\end{theorem}

The positive assertion follows from the local smoothing theorem of Guth,
H. Wang, and R. Zhang and the transfer principle; it is also stated for generalized
spherical means in \cite{LSSY2024}.  The necessary condition is the
specialization of Theorem~\ref{thm:LSSY} to $n=2$.  For real orders, the
radial results in \cite{NRS2023} already imply the same obstruction, since
boundedness on all functions entails boundedness on the radial subspace. Equivalently, define
$$
 E_2(p):=
 \min\Big\{\frac12-\frac1p,\frac1p\Big\}
 =
 \begin{cases}
  \dfrac12-\dfrac1p, & 2<p\le4,\\
  \dfrac1p, & 4\le p<\infty.
 \end{cases}
$$
Then the positive result is $\operatorname{Re}\alpha>-E_2(p)$, while
necessity gives $\operatorname{Re}\alpha\ge-E_2(p)$.  In the
$(1/p,\operatorname{Re}\alpha)$-plane, the two boundary segments meet at
the corner
$$
 \Big(\frac1p,\operatorname{Re}\alpha\Big)
 =\Big(\frac14,-\frac14\Big).
$$

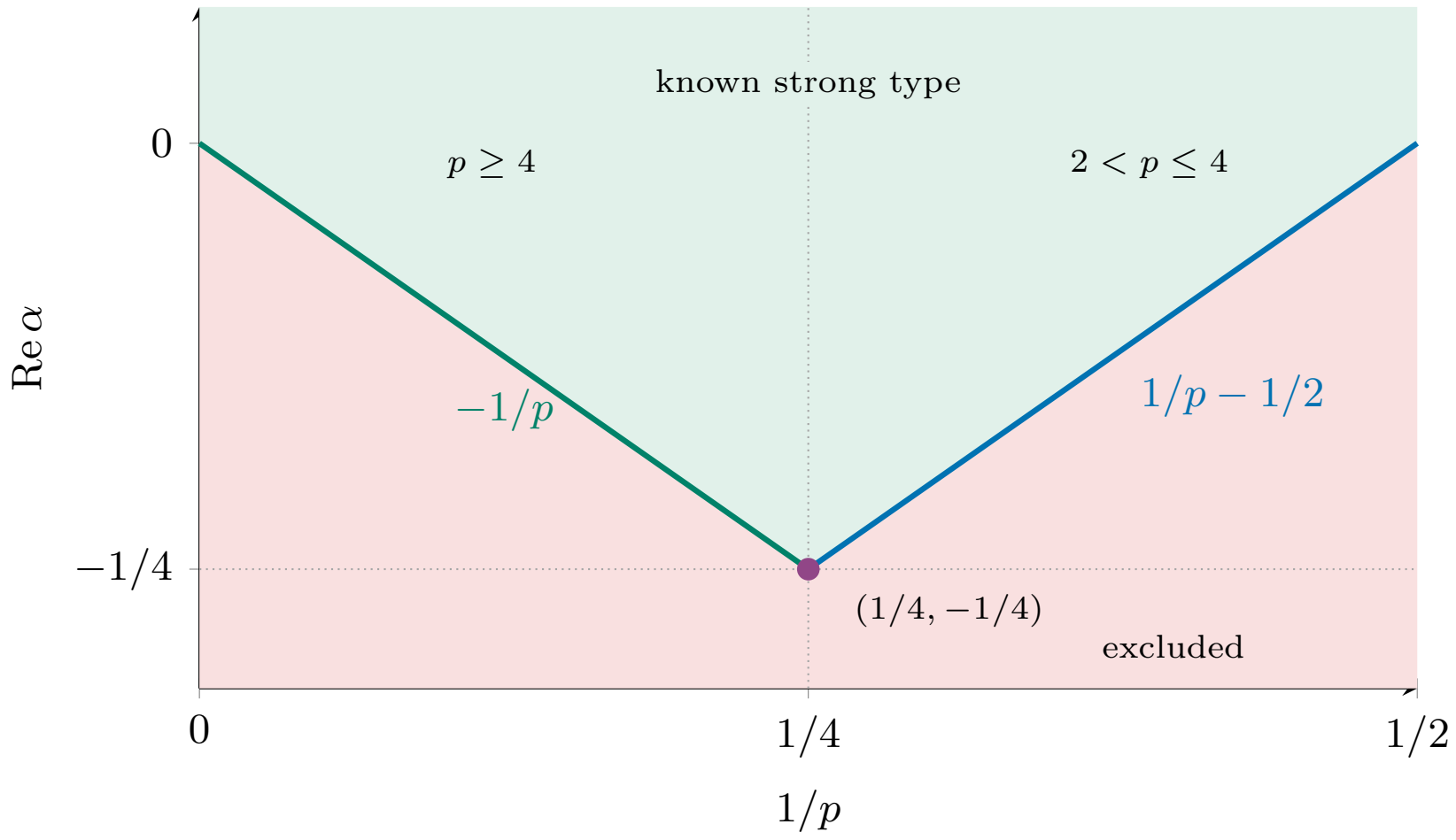
\begin{figure}[htbp]
\centering
\begin{tikzpicture}
\begin{axis}[
  width=0.64\textwidth,  
  height=0.40\textwidth, 
  xmin=0,
  xmax=0.5,
  ymin=-0.32,
  ymax=0.08,
  axis lines=left,
  xlabel={$1/p$},
  ylabel={$\operatorname{Re}\alpha$},
  xtick={0,0.25,0.5},
  xticklabels={$0$,$1/4$,$1/2$},
  ytick={-0.25,0},
  yticklabels={$-1/4$,$0$},
  tick label style={font=\small},
  label style={font=\small},
  clip=false
]
\path[draw=none,fill=KnownFill]
  (axis cs:0,0.08) -- (axis cs:0.5,0.08) -- (axis cs:0.5,0)
  -- (axis cs:0.25,-0.25) -- (axis cs:0,0) -- cycle;
\path[draw=none,fill=ExcludedFill]
  (axis cs:0,-0.32) -- (axis cs:0,0) -- (axis cs:0.25,-0.25)
  -- (axis cs:0.5,0) -- (axis cs:0.5,-0.32) -- cycle;

\addplot[forget plot,black!35,thin,densely dotted]
  coordinates {(0.25,-0.32) (0.25,0.08)};
\addplot[forget plot,black!35,thin,densely dotted]
  coordinates {(0,-0.25) (0.5,-0.25)};

\addplot[GHLWGreen,very thick]
  coordinates {(0,0) (0.25,-0.25)}
  node[pos=0.50,below=3pt,font=\small] {$-1/p$};
\addplot[SteinBlue,very thick]
  coordinates {(0.25,-0.25) (0.5,0)}
  node[pos=0.50,anchor=north west,xshift=2pt,font=\small] {$1/p-1/2$}; 
\addplot[only marks,mark=*,mark size=2.2pt,BoundaryPurple]
  coordinates {(0.25,-0.25)};

\node[font=\scriptsize,fill=KnownFill,inner sep=1.5pt]
  at (axis cs:0.25,0.035) {known strong type};
\node[font=\scriptsize,fill=ExcludedFill,inner sep=1.5pt]
  at (axis cs:0.40,-0.295) {excluded}; 
\node[font=\scriptsize,anchor=north west]
  at (axis cs:0.262,-0.255) {$(1/4,-1/4)$}; 
\node[font=\scriptsize,anchor=south]
  at (axis cs:0.12,-0.03) {$p\ge4$};
\node[font=\scriptsize,anchor=south]
  at (axis cs:0.39,-0.03) {$2<p\le4$};
\end{axis}
\end{tikzpicture}
\caption{The complete open boundedness region in dimension two.  Above the
two-segment boundary, $M^\alpha$ is bounded on $L^p(\R^2)$; below it,
boundedness is impossible.  The corner at $(1/4,-1/4)$ corresponds to $p=4$.}
\label{fig:planar-phase}
\end{figure}

This is shown in Figure~\ref{fig:planar-phase}: the blue segment $2<p<4$ is controlled by
interpolation with the $L^2$ estimate, whereas the green segment $p\ge4$ comes from
the sharp planar local smoothing theorem.

\subsection{Comparison with the result of Miao, Yang, and Zheng}

Specializing \eqref{eq:MYZepsilon} to $n=2$ gives
$$
 E_{\rm MYZ}(p,2)
 =\min\Big\{\frac14-\frac1{2p},\frac1p\Big\}.
$$
Its two branches meet at $p=6$.  Thus the former sufficient condition can be
written
$$
 \operatorname{Re}\alpha>
 \begin{cases}
  -\dfrac14+\dfrac1{2p}, & 2<p\le6,\\
  -\dfrac1p, & 6\le p<\infty.
 \end{cases}
$$
The gain obtained by the planar theorem can be seen in the
$E$-notation:
$$
 E_2(p)-E_{\rm MYZ}(p,2)
 =
 \begin{cases}
  \dfrac14-\dfrac1{2p}, & 2<p\le4,\\
  \dfrac3{2p}-\dfrac14, & 4\le p<6,\\
  0, & 6\le p<\infty.
 \end{cases}
$$
On $4\le p<6$, this is the direct gain from the sharp cone square function
estimate.  On $2<p<4$, it is an indirect gain obtained by interpolating the
new endpoint $p=4$ with the energy estimate at $p=2$. As a consistency check, setting $\alpha=0$ in
Theorem~\ref{thm:n2-open-region} gives boundedness of the circular maximal
operator for every $p>2$, recovering Bourgain's theorem
\cite{Bourgain1986}.

\subsection{What the word ``open'' excludes}

Theorem~\ref{thm:n2-open-region} classifies every point away from the two
boundary segments: strict inequality above the boundary gives strong type,
whereas strict inequality below it is impossible.  As explained in
Remark~\ref{rem:strict-inequalities}, the transfer argument makes no
assertion at equality. The radial theory illustrates why the boundary should be kept separate.
For real $\alpha$ and $2<p\le4$,
\cite[Corollary~1.6]{NRS2023} shows at
$$
 \alpha=\frac1p-\frac12
$$
restricted weak type but failure of weak type, and hence failure of strong
type.  For $p>4$, by contrast, the radial operator is of strong type at
$\alpha=-1/p$ \cite[Corollary~1.5]{NRS2023}.  Thus radial functions provide
no obstruction on the second open boundary segment; see
Section~\ref{sec:open}.

\section{\texorpdfstring{The local smoothing theorems of Gan, He, Li, and Wu and of X. Wang}{The local smoothing theorems of Gan, He, Li, and Wu and of X. Wang}}
\label{sec:GHLW}

The recent theorems of Gan, He, Li, and Wu \cite{GHLW2026} and of
X.~Wang \cite{Wang2026} are most naturally
stated in the FIO framework of Section~\ref{sec:dictionary}.   The Euclidean local smoothing conjecture begins at
$$
 q_n=2+\frac2{n-1}.
$$
For general compact manifolds, the examples of Minicozzi and Sogge
\cite{MinicozziSogge1997} give the larger necessary threshold
\begin{equation}
\label{eq:pn-plus}
 p_n^+:=
 \begin{cases}
  2+\dfrac8{3n-3}, & n\text{ odd},\\[6pt]
  2+\dfrac8{3n-2}, & n\text{ even}.
 \end{cases}
\end{equation}
Finally, the former threshold supplied by Bourgain and Demeter is
$$
 p_{\rm BD}(n)=2+\frac4{n-1}.
$$
The Euclidean threshold $q_n$ concerns the constant-coefficient half-wave,
whereas $p_n^+$ is an obstruction for the variable-coefficient problem.

The best threshold currently available is
\begin{equation}
\label{eq:Pn}
 P_n:=
 \begin{cases}
  2+\dfrac8{3n-3}, & n\text{ odd},\\[7pt]
  2+\dfrac8{3n-2-\omega(n)}, & n\text{ even},
 \end{cases}
 \qquad
 \omega(n):=\frac{12(n+2)}{2n^2+7n+14}.
\end{equation}
The odd case is proven in \cite{GHLW2026}, where in even dimensions the authors
obtained $2+8/(3n-2-\omega_0(n))$ with
$\omega_0(n)=4(7n-2)/(3n^2+5n+14)$; the smaller value in \eqref{eq:Pn} is
presented in \cite{Wang2026}.
Thus $P_n=p_n^+$ when $n$ is odd.  When $n$ is even,
$P_n>p_n^+$, but
$$
 P_n=p_n^++O(n^{-3}).
$$
For every $n\ge3$,
\begin{equation}
\label{eq:threshold-ordering}
 q_n<p_n^+\le P_n<p_{\rm BD}(n),
\end{equation}
where equality in the middle occurs precisely in odd dimensions.  The first
few values are
$$
 P_3=\frac{10}{3},
 \qquad
 P_4=\frac{482}{167},
 \qquad
 P_5=\frac83.
$$

\subsection{The theorem in the form needed here}

Let $\mathcal F$ be the order-zero FIO in \eqref{eq:FIO-model}, with an
amplitude compactly supported in $(x,t)$ and a nondegenerate phase satisfying
the positive-definite curvature condition described there.  In even
dimensions we assume in addition the quadratic-cone condition of
\cite{Wang2026}: for each $z$ there is a nondegenerate symmetric matrix
$Q_z$, independent of $\xi$, such that
$\partial_z\phi(z;\xi)^{T}Q_z\,\partial_z\phi(z;\xi)=0$.  By the eikonal
equation (cf. \cite{Wang2026}) this holds for Riemannian half-wave parametrices, and in particular
for the Euclidean phase $x\cdot\xi+t|\xi|$ with $Q=\operatorname{diag}(I,-1)$.

\begin{theorem}[Gan, He, Li, and Wu \cite{GHLW2026}; X.~Wang \cite{Wang2026}]
\label{thm:GHLW}
Let $n\ge3$ and $P_n\le p<\infty$.  For every
\begin{equation}
\label{eq:GHLW-regularity}
 s>\gamma(p,n)-\frac1p
 =\frac{n-1}{2}-\frac np,
\end{equation}
one has
$$
 \|\mathcal Ff\|_{L^p(\R^n\times[1,2])}
 \lesssim_s
 \|f\|_{W^{s,p}(\R^n)}.
$$
In particular,
$$
 \big\|e^{\pm it|D|}f\big\|_{L^p(\R^n\times[1,2])}
 \lesssim_s
 \|f\|_{W^{s,p}(\R^n)}
$$
under the same assumptions.
\end{theorem}

For the
constant-coefficient half-wave, translation invariance and a partition of
unity give the global spatial estimate on the fixed interval $[1,2]$.   The FIO
formulation also permits uniformly bounded order-zero amplitudes, which is the
stability required in \eqref{eq:transfer-LS-piece} after the Bessel
In odd dimensions, the threshold for Theorem~\ref{thm:GHLW} is
optimal for general compact manifolds because $P_n=p_n^+$.  In even
dimensions only the small gap between $p_n^+$ and $P_n$ remains for the
variable-coefficient problem; the quadratic-cone condition costs
nothing here, since the Euclidean phase satisfies it, and only for general
FIOs without it is the threshold of \cite{GHLW2026} with $\omega_0(n)$ the
best known.  For the Euclidean half-wave, however, the
conjectural starting point is still the smaller exponent $q_n$.  Thus the new
results improve the former threshold from $p_{\rm BD}(n)$ to $P_n$, but do
not settle Euclidean local smoothing on $q_n\le p<P_n$ yet.

\subsection{\texorpdfstring{Main ideas of the proof of Theorem~\ref{thm:GHLW}}{Main ideas of the proof}}

Several
features of the proof in \cite{GHLW2026} explain why the threshold improves and why the answer
depends on the parity of the dimension.

\begin{enumerate}
\item \emph{A scale-invariant local estimate.}
After frequency localization and rescaling, the problem is reduced
schematically to
$$
 \|\mathcal Ff\|_{L^p(B_R^{n+1})}
 \lesssim_\eps
 R^{\gamma(p,n)+\eps}\|f\|_{L^p},
$$
when $\widehat f$ is supported in a fixed annulus.  A wave packet decomposition writes
$f=\sum_{T\in\mathbb T}f_T$.  The corresponding space--time packets are
concentrated on conic planks of dimensions
$$
 1\times R^{1/2}\times\cdots\times R^{1/2}\times R.
$$

\item \emph{Wave packet density.}
A direct induction with the $L^p$ norm of the data is not compatible with
Lorentz rescaling.  The authors replace it by a mixed quantity of the form
$$
 \|f\|_2^{2/p}
 \mathcal W(f;B_R^{n+1})^{1-2/p}.
$$
Here $\mathcal W$ is the \emph{wave packet density} introduced in
\cite[Definition~4.1]{GHLW2026}. For a ball centered at the origin, it
is\footnote{We follow the notation of
\cite[Definition~4.1]{GHLW2026}, except that the scale is denoted
$\sigma$ instead of $s$, which is reserved here for Sobolev regularity,
and that the dependence on the center of the ball is suppressed.}
\begin{equation*}
 \mathcal W(f;B_R^{n+1})
 :=\sup_{R^{-1/2}\le\sigma\le1}\sup_{\tau\in\Theta_\sigma}
   \sup_{V\parallel V_{\tau,R}}
 \bigg(\frac1{|V|}
 \sum_{\substack{\theta\subset\tau\\\theta\in\Theta_{R^{-1/2}}}}
 \sum_{\substack{T\in\mathbb T_\theta\\T^\flat\subset V}}
 \|f_T\|_2^2\bigg)^{1/2}.
\end{equation*}
In this formula, $\Theta_\sigma$ is a covering by angular caps of aperture
$\sigma$, $\mathbb T_\theta$ consists of packets with direction in
$\theta$, and $T^\flat$ is the spatial base of the packet in the local
coordinates used for the decomposition. The box $V_{\tau,R}$ has dimensions
$R\sigma^2\times R\sigma\times\cdots\times R\sigma$ and orientation
determined by $\tau$; $V\parallel V_{\tau,R}$ ranges over its translates.
Thus $\mathcal W$ measures the largest square root of the packet energy
per unit spatial volume across all these scales and boxes. It replaces
$L^\infty$ control in a way compatible with orthogonality and Lorentz rescaling,
which makes induction on scales possible.

\item \emph{Broadness and two ends.}
The broad--narrow decomposition separates points receiving contributions
from several transverse packet directions from those where the directions
cluster.  Narrow contributions are treated recursively after Lorentz
rescaling.  In the broad part, a two-ends condition prevents the portion of a
packet under consideration from concentrating near a single small region.
Since the two conditions cannot be imposed independently by a single
pigeonholing argument, they are obtained by an iterative selection
procedure.

\item \emph{Refined decoupling, orthogonality, and incidence geometry.}
Refined cone decoupling controls the broad contribution at the cone
decoupling exponent
$$
 p_{\rm dec}(n):=\frac{2(n+1)}{n-1}=p_{\rm BD}(n).
$$
We use the notation $p_{\rm dec}(n)$ here because the paper denotes this
quantity by $q_n$, while in this survey $q_n=2n/(n-1)$ is reserved for the
Euclidean local smoothing threshold.  The refined decoupling estimate is
combined with $L^2$ orthogonality and with incidence bounds between unit
balls and conic planks.  The hairbrush geometry of the cone controls the
multiplicity of packets through a ball once the two-ends and broadness
conditions are available.

\item \emph{Dimension-sensitive closure of the induction.}
The broad norm is chosen with
$$
 k_n=
 \begin{cases}
  (n+5)/2, & n\text{ odd},\\
  (n+4)/2, & n\text{ even}.
 \end{cases}
$$
In dimension $n=3$, a covering argument gives a refined $L^2$ estimate for
packets inside slabs.  In dimension $n=5$, a trilinear estimate exploits the
transversality of the relevant slab normals.  In odd dimensions $n\ge7$,
the required exponent follows by interpolating the refined decoupling,
$L^2$, and $k_n$-broad estimates.  The even-dimensional argument does not reach $p_n^+$; its loss is encoded by $\omega_0(n)$. X. Wang
\cite{Wang2026} reduces it to $\omega(n)$ in \eqref{eq:Pn} by separating the
broad structure used for the $L^p$ bound from the one used in the hairbrush
argument, and by adding a multilinear estimate.
\end{enumerate}

For the maximal application, the outcome of this strategy is the input \eqref{eq:GHLW-regularity}.  Below $P_n$, it is combined with
the $L^2$ estimate by interpolation before the transfer principle is applied.

\section{A new maximal corollary and the interpolation below
\texorpdfstring{$P_n$}{Pn}}
\label{sec:new-corollary}

To use Theorem~\ref{thm:GHLW} below $P_n$, interpolate its estimate at
$p=P_n$ with the $L^2$ identity.  If $2<p<P_n$, write
$$
 \frac1p=\frac{1-\theta}{2}+\frac{\theta}{P_n}.
$$
Since the regularity at $P_n$ may be taken arbitrarily close to
$(n-1)/2-n/P_n$, the interpolated regularity may be taken arbitrarily close
to
$$
 \theta\Big(\frac{n-1}{2}-\frac n{P_n}\Big)
 =
 c_n\Big(\frac12-\frac1p\Big),
$$
where
\begin{equation}
\label{eq:cn}
 c_n:=\frac{(n-1)P_n-2n}{P_n-2}.
\end{equation}
The inequality $P_n<p_{\rm BD}(n)$ implies
$$
 c_n<\frac{n-1}{2}<n.
$$
Thus, for $2<p<P_n$, one has local smoothing for every
\begin{equation}
\label{eq:sinterp}
 s>s_{\rm int}(p,n)
 :=c_n\Big(\frac12-\frac1p\Big).
\end{equation}
The limiting value at $p=P_n$ agrees with
$(n-1)/2-n/P_n$.

\begin{corollary}[Theorem~\ref{thm:GHLW} and the transfer principle]
\label{cor:new-maximal}
Let $n\ge3$ and $p>2$.  Then $M^\alpha$ is bounded on $L^p(\R^n)$ whenever
$$
 \operatorname{Re}\alpha>-E_{\rm new}(p,n),
$$
where
\begin{equation}
\label{eq:Enew}
 E_{\rm new}(p,n)=
 \begin{cases}
 \displaystyle
 \frac{n-1}{2}-\frac1p
 -c_n\Big(\frac12-\frac1p\Big),
 &2<p<P_n,\\[9pt]
 \displaystyle\frac{n-1}{p},&p\ge P_n.
 \end{cases}
\end{equation}
Here $P_n$ and $c_n$ are given by \eqref{eq:Pn} and \eqref{eq:cn}.
\end{corollary}

\begin{proof}
For $p\ge P_n$, insert \eqref{eq:GHLW-regularity} into
Proposition~\ref{prop:transfer}.  The local smoothing and fixed-time
contributions then agree:
$$
 \frac{n-1}{2}-\frac1p
 -\Big(\frac{n-1}{2}-\frac np\Big)
 =\frac{n-1}{p}.
$$
This gives the second line of \eqref{eq:Enew}.

Let $2<p<P_n$.  Inserting \eqref{eq:sinterp} into the transfer principle
gives the minimum of
$$
 \frac{n-1}{p}
 \qquad\text{and}\qquad
 \frac{n-1}{2}-\frac1p
 -c_n\Big(\frac12-\frac1p\Big).
$$
The first quantity minus the second equals
$$
 (n-c_n)\Big(\frac1p-\frac1{P_n}\Big)>0.
$$
Hence the second quantity is the relevant one, which proves the first line
of \eqref{eq:Enew}.  The strict inequality in
$\operatorname{Re}\alpha$ absorbs the arbitrarily small loss in the local
smoothing estimate.
\end{proof}

Table~\ref{tab:small-dimensional-thresholds} records the first three dimensional cases.
Besides the new threshold $P_n$, it includes the Euclidean conjectural
threshold $q_n$, the manifold obstruction $p_n^+$, and the former decoupling
threshold $p_{\rm BD}(n)$.  The equality $P_n=p_n^+$ in odd dimensions and
the small even-dimensional discrepancy are visible in these values.

\begin{table}[htbp]
\centering
\caption{Thresholds and the new interpolated depth in small dimensions.}
\label{tab:small-dimensional-thresholds}
\small
\setlength{\tabcolsep}{6pt}
\renewcommand{\arraystretch}{1.35}
\begin{tabular}{cccccc}
\toprule
$n$ & $q_n$ & $p_n^+$ & $P_n$ & $p_{\rm BD}(n)$
& $E_{\rm new}(p,n)$, $2<p<P_n$\\
\midrule
$3$ & $3$ & $10/3$ & $10/3$ & $4$
& $\displaystyle \frac34-\frac1{2p}$\\
$4$ & $8/3$ & $14/5$ & $482/167$ & $10/3$
& $\displaystyle \frac{167}{148}-\frac{19}{74p}$\\
$5$ & $5/2$ & $8/3$ & $8/3$ & $3$
& $\displaystyle \frac32$\\
\bottomrule
\end{tabular}
\end{table}

The same inequality also quantifies the improvement over
Theorem~\ref{thm:MYZ}: for $2<p<P_n$,
\begin{equation}
\label{eq:gain}
 E_{\rm new}(p,n)-E_{\rm MYZ}(p,n)
 =\Big(\frac{n-1}{2}-c_n\Big)
   \Big(\frac12-\frac1p\Big)>0.
\end{equation}
For $P_n\le p<p_{\rm BD}(n)$ the new results move directly to the sharp
branch $(n-1)/p$.

\begin{figure}[htbp]
\centering
\begin{tikzpicture}
\begin{axis}[
  width=0.72\textwidth,  
  height=0.42\textwidth, 
  xmin=0.23,
  xmax=0.5,
  ymin=-0.69,
  ymax=-0.43,
  axis lines=left,
  xlabel={$1/p$},
  ylabel={$\operatorname{Re}\alpha$},
  xtick={0.25,0.3,0.333333,0.5},
  xticklabels={
    $1/p_{\rm BD}(3)$,
    $1/P_3$,
    $1/q_3$,
    $1/2$
  },
  x tick label style={font=\scriptsize,rotate=32,anchor=north east},
  ytick={-0.666667,-0.6,-0.5},
  yticklabels={$-2/3$,$-3/5$,$-1/2$},
  y tick label style={font=\small},
  label style={font=\small},
  legend style={
    at={(0.02,0.98)},
    anchor=north west,
    draw=none,
    fill=white,
    fill opacity=0.92,
    text opacity=1,
    font=\scriptsize,
    legend columns=3,
    column sep=7pt
  },
  clip=false
]
\path[draw=none,fill=GainFill]
  (axis cs:0.25,-0.5) -- (axis cs:0.3,-0.6)
  -- (axis cs:0.5,-0.5) -- cycle;
\path[draw=none,fill=GapFill]
  (axis cs:0.3,-0.6) -- (axis cs:0.5,-0.5)
  -- (axis cs:0.333333,-0.666667) -- cycle;

\addplot[forget plot,black!35,thin,densely dotted]
  coordinates {(0.25,-0.69) (0.25,-0.43)};
\addplot[forget plot,black!35,thin,densely dotted]
  coordinates {(0.3,-0.69) (0.3,-0.43)};
\addplot[forget plot,black!35,thin,densely dotted]
  coordinates {(0.333333,-0.69) (0.333333,-0.43)};

\addplot[MYZOrange,very thick]
  coordinates {(0.23,-0.46) (0.25,-0.5) (0.5,-0.5)};
\addlegendentry{Miao, Yang, and Zheng}

\addplot[GHLWGreen,very thick]
  coordinates {(0.23,-0.46) (0.3,-0.6) (0.5,-0.5)};
\addlegendentry{New boundary}

\addplot[BoundaryPurple,very thick]
  coordinates {(0.23,-0.46) (0.333333,-0.666667) (0.5,-0.5)};
\addlegendentry{Necessary boundary}

\node[font=\scriptsize,fill=GainFill,inner sep=1.5pt]
  at (axis cs:0.276,-0.525) {new gain};
\node[font=\scriptsize,fill=GapFill,inner sep=1.5pt]
  at (axis cs:0.405,-0.585) {remaining gap};
\end{axis}
\end{tikzpicture}
\caption{A zoom of the higher-dimensional diagram for $n=3$.  The blue
region is the improvement over the Miao, Yang, and Zheng theorem; the yellow
region is the gap that remains for $2<p<P_3$.  At $p=P_3$ the new boundary
meets the necessary line, while $p_{\rm BD}(3)$ marks the former point of
contact.  The exponent $q_3$ is the conjectured starting point for sharp
Euclidean local smoothing.}
\label{fig:higher-dimensional-zoom}
\end{figure}
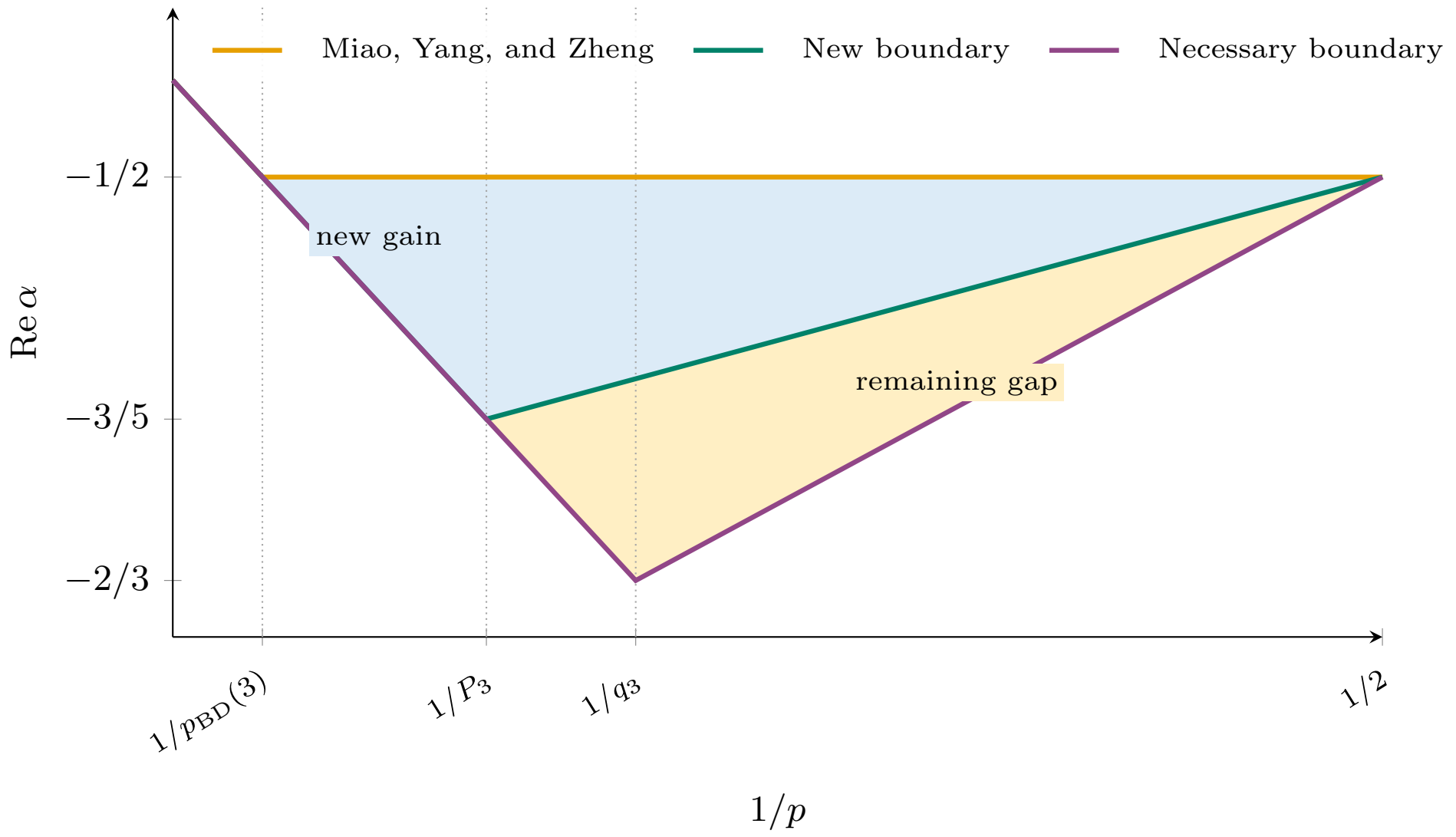

Figure~\ref{fig:higher-dimensional-zoom} enlarges the part of the diagram in
which the improvement occurs.  We again use $n=3$, so that all three
thresholds have simple values.  Reading the horizontal axis from right to
left corresponds to increasing $p$: $q_3=3$ is the conjectural Euclidean
threshold, $P_3=10/3$ is the new local smoothing threshold, and
$p_{\rm BD}(3)=4$ is the former one.

Consequently the general maximal problem is now sharp, apart from boundary
questions, for every $p\ge P_n$.  The sufficient and necessary regions can
still differ only when
$$
 2<p<P_n.
$$

The earlier higher-dimensional local smoothing improvement of Gao,
Liu, Miao, and Xi \cite{GLMX2023} forms part of the development between
cone decoupling and the estimates discussed above. We do not need its
separate exponent here, since the later results of Gan, He, Li, and Wu
and of X. Wang yield the stronger maximal region shown in the diagram.

\section{Further consequences and related directions}
\label{sec:further}

The analytic family has several useful interpretations that are not needed
for the proof of the phase diagram but help place the problem in a broader
context.  We record three of them: derivatives of ordinary spherical means,
pointwise convergence, and local improving estimates.

\subsection{Derivatives of spherical means}

Let
$$
 A_t f(x):=\int_{\mathbb S^{n-1}}f(x+ty)\,d\sigma(y)=\M_t^0f(x),
$$
where $d\sigma$ is normalized surface measure, and for
$k=0,1,2,\ldots$ define
$$
 A_t^{(k)}f(x):=t^k\Big(\frac d{dt}\Big)^kA_tf(x),
 \qquad
 A_*^{(k)}f(x):=\sup_{t>0}|A_t^{(k)}f(x)|.
$$
Kwa\'snicki \cite{Kwasnicki2026} observed that the operators
$A_t^{(k)}$ can be recovered from the generalized means of integer negative
order.  More precisely, whenever none of the orders $-j$, $0\le j\le k$,
belongs to the exceptional set $\{-n/2,-n/2-1,\ldots\}$ (see
Section~\ref{sec:intro}), that is, when $n$ is odd or $n>2k$,
there are constants $c_{n,k,j}$, with $c_{n,k,k}\ne0$, such that
\begin{equation}
\label{eq:derivative-triangular}
 A_t^{(k)}f=\sum_{j=0}^k c_{n,k,j}\M_t^{-j}f.
\end{equation}
Since $c_{n,k,k}\ne0$, the relation can be inverted: $\M_t^{-k}f$ is
in turn a linear combination of $A_t^{(0)}f,\ldots,A_t^{(k)}f$. At the
level of Bessel multipliers, this is an iteration of
$$
 \frac d{dr}\big[r^\nu J_\nu(r)\big]=r^\nu J_{\nu-1}(r).
$$
Indeed, writing $r^{-\nu}J_\nu(r)=r^{-2\nu}\,r^\nu J_\nu(r)$ and using
$c_{n,\alpha}=2\nu\,c_{n,\alpha-1}$, with $\nu=n/2+\alpha-1$ as in
Section~\ref{sec:dictionary}, one gets
$$
 t\partial_t\M_t^\alpha f=(n-2+2\alpha)\big(\M_t^{\alpha-1}f-\M_t^\alpha f\big);
$$
for $\alpha=0$ and $n\ge3$ this reads $A_t^{(1)}f=(n-2)(\M_t^{-1}f-A_tf)$.
Thus decreasing the order of Stein's family by one corresponds, up to
lower-order terms, to applying one scale-invariant radial derivative to the
ordinary spherical mean.  Formula \eqref{eq:derivative-triangular} makes the
negative values of $\alpha$ concrete by measuring differentiability in the
radius.

The resulting maximal theorem can be stated without reference to analytic
continuation.

\begin{theorem}[Kwa\'snicki \cite{Kwasnicki2026}]
\label{thm:Kwasnicki}
Let $k\ge0$, $n\ge2k+3$, and
$$
 \frac{n}{n-k-1}<p<\frac{n-2}{k},
$$
where $(n-2)/k=\infty$ if $k=0$.  Then
$$
 \|A_*^{(k)}f\|_{L^p(\R^n)}\le C_{k,p}\|f\|_{L^p(\R^n)},
$$
where the constant is independent of $n$.
\end{theorem}

The two restrictions on $p$ are Stein's conditions
\eqref{eq:Stein-region} at $\alpha=-k$: the lower endpoint is the
positive-kernel obstruction for $p<2$, and the upper one is the line
$\operatorname{Re}\alpha=(2-n)/p$.
The dimension-free constant is obtained by expressing a spherical mean in
dimension $n$ as an average of means on rotated copies of
$\R^{n-1}$.  This averaging identity commutes with radial differentiation,
so a bound in one dimension propagates to every higher dimension with the
same constant.  By Theorem~\ref{thm:LSSY}, no $p>(n-1)/k$ is
admissible, so the conjectured range is $n/(n-k-1)<p<(n-1)/k$.  It is
obtained in \cite{Kwasnicki2026} for $n\ge2k+3$, with constants depending
on $n$, and the case $n=2k+2$, where the range is $2<p<2+1/k$, is listed
there as the only open one.  Corollary~\ref{cor:new-maximal} settles the
first instance.

\begin{corollary}
\label{cor:n4}
The maximal operator $A_*^{(1)}$ is bounded on $L^p(\R^4)$ for
$2<p<3$.
\end{corollary}

\begin{proof}
By Table~\ref{tab:small-dimensional-thresholds},
$E_{\rm new}(p,4)=167/148-19/(74p)>1$ for $2<p<P_4$, and
$E_{\rm new}(p,4)=3/p>1$ for $P_4\le p<3$.  Hence $M^{-1}$ is bounded on
$L^p(\R^4)$ for $2<p<3$ by Corollary~\ref{cor:new-maximal}, and so is
$A_*=M^0$.  The claim follows from
$A_t^{(1)}f=2(\M_t^{-1}f-A_tf)$.
\end{proof}

The range is optimal except for the endpoint $p=3$: at $p=2$ the
operator $M^{-1}$ is unbounded already on radial functions
\cite[Corollary~1.6]{NRS2023}, and the identity above transfers this to
$A_*^{(1)}$.  The argument does not extend to $k\ge2$, and the
obstruction can be made precise.  Let $n=2k+2$.  Since $-k=1-n/2$, the
point $(1/p,-k)$ lies on the conjectured boundary at $p=2$, and
Proposition~\ref{prop:transfer} covers it for $2<p<2+1/k$ exactly when
local smoothing holds with some $s<1/2-1/p$.  By interpolation with
$L^2$, this is the case if and only if the sharp estimate
\eqref{eq:GHLW-regularity} is available at some exponent
$p_0<2+1/k=2+2/(n-2)$, that is, if and only if $c_n<1$.  For $n=4$ this
holds because $P_4<3$.  For even $n\ge6$ one has $2+2/(n-2)\le p_n^+$,
so the required input lies below the Minicozzi--Sogge threshold: it is a
Euclidean local smoothing estimate, which no variable-coefficient
method can provide, and which is what is missing in
Section~\ref{sec:open}.  Kwa\'snicki's dimensional induction does not help
either, since it only passes from $\R^{n-1}$ to $\R^n$.  Thus, for
$k\ge2$, the cases $n=2k+2$ are equivalent to progress on the Euclidean
conjecture below $p_n^+$; they hold on radial functions by
Theorem~\ref{thm:NRS-radial}.

There is also a fractional version.  For $\operatorname{Re}\alpha>0$,
integration in polar coordinates gives
$$
 \M_t^\alpha f(x)=\frac{2\Gamma(\alpha+n/2)}{\Gamma(\alpha)\Gamma(n/2)}
 \int_0^1(1-r^2)^{\alpha-1}r^{n-1}A_{rt}f(x)\,dr,
$$
so that, up to the constant $\Gamma(\alpha+n/2)/\Gamma(n/2)$,
$\M_t^\alpha f(x)$ is the Erd\'elyi--Kober fractional integral of order
$\alpha$, with respect to $t^2$, of $t\mapsto A_tf(x)$; see
\cite{Sneddon1975} and \cite[Ch.~1]{Kiryakova1994}.  For negative orders,
$\M_t^\alpha f$ is the analytic continuation of this integral, that is, an
Erd\'elyi--Kober fractional derivative of order $-\alpha$.  Integer
derivatives are therefore the discrete part of a fractional calculus in
the dilation parameter.

\subsection{Maximal bounds and pointwise convergence}

Our normalization gives $\widehat m_\alpha(0)=1$.  Hence, for Schwartz
functions and every admissible order,
$$
 \lim_{t\to0}\M_t^\alpha f(x)=f(x).
$$
When $\operatorname{Re}\alpha\le0$, this assertion follows from the Bessel multiplier formula and its
expansion at the origin. A strong $L^p$ estimate for $M^\alpha$, together with density of the Schwartz
class, upgrades this convergence to
\begin{equation}
\label{eq:ae-convergence}
 \lim_{t\to0}\M_t^\alpha f(x)=f(x)
 \quad\text{for almost every }x,
 \qquad f\in L^p(\R^n).
\end{equation}
Thus each interior point of the phase
diagram also gives a differentiation theorem for a possibly distributional
family of means.

\subsection{Local \texorpdfstring{$L^p\to L^q$}{Lp to Lq} estimates}

Off-diagonal estimates become possible once the radii are restricted to a
compact interval.  Set
$$
 M_{[1,2]}^\alpha f(x):=\sup_{1\le t\le2}|\M_t^\alpha f(x)|.
$$
For the global operator $M^\alpha$, dilation invariance forces $p=q$ in any
unweighted estimate $L^p\to L^q$.  The restriction $1\le t\le2$ removes this
scaling obstruction and allows genuine improving estimates.

For $1\le p,q<\infty$, Liu, Shen, Song, and Yan
\cite{LSSYimproving2025} introduced the function
\begin{equation}
\label{eq:sigma-pq}
 \sigma_n(p,q):=\max\Big\{
 \frac1p-\frac n q,
 \frac{n+1}{2p}-\frac{n-1}{2}\Big(\frac1q+1\Big),
 \frac n p-n+1
 \Big\}.
\end{equation}
Their result says that an estimate
$$
 \|M_{[1,2]}^\alpha f\|_{L^q(\R^n)}\lesssim
 \|f\|_{L^p(\R^n)}
$$
requires $q\ge p$ and
$\operatorname{Re}\alpha\ge\sigma_n(p,q)$.  Conversely, for $q\ge p$ it
holds in dimension two when
$$
 \operatorname{Re}\alpha>\sigma_2(p,q),
$$
and in dimensions $n>2$ when
\begin{equation}
\label{eq:LSSY-improving-sufficient}
 \operatorname{Re}\alpha>
 \max\Big\{\sigma_n(p,q),
 \frac1{2p}-\frac{n-2}{2q}-\frac{n-1}{4}\Big\}.
\end{equation}
The three terms in \eqref{eq:sigma-pq} arise from different test functions
and geometries.  On the diagonal $p=q$, they reduce respectively to
$$
 -\frac{n-1}{p},
 \qquad
 \frac1p-\frac{n-1}{2},
 \qquad
 1-n+\frac np.
$$

A complementary earlier contribution is due to Bagchi, Hait, and Senthil
Raani \cite{BHR2022}.  They established $L^p\to L^q$ improving estimates for
unit-scale generalized spherical averages, together with the quantitative
translation-continuity estimates needed in Lacey's sparse domination method \cite{Lacey2019}.
Using these two inputs, they obtained sparse bounds for the associated
lacunary and full generalized spherical maximal operators, and consequently
weighted as well as unweighted estimates.  

The local operator has also been studied in the scale of Hardy spaces
for Fourier integral operators $\mathcal H^{s,p}_{FIO}(\R^n)$.  Ghosh, Liu,
Rozendaal, and Song \cite{GLRS2026} prove
$\|M_{[1,2]}^\alpha f\|_{L^p}\lesssim\|f\|_{\mathcal H^{s,p}_{FIO}}$ with a
regularity index $s$ that is sharp for
$p\in[1,2]\cup[p_{\rm BD}(n),\infty)$ and which, through the Sobolev
embeddings for these spaces, recovers for $M_{[1,2]}^\alpha$ the $L^p$
regions of Theorems~\ref{thm:Stein} and~\ref{thm:MYZ}; for $\alpha=0$ the
method also covers geodesic spheres on compact manifolds.

One may restrict the radii further.  For $E\subset[1,2]$, define
$$
 M_E^\alpha f(x):=\sup_{t\in E}|\M_t^\alpha f(x)|.
$$
Already at $\alpha=0$, this family interpolates between a fixed spherical
average when $E$ is a singleton and the local spherical maximal operator when
$E=[1,2]$.  For fractal sets of radii, the sharp improving region depends not
only on the upper Minkowski dimension of $E$ but also on finer concentration
data encoded by Assouad-type spectra; see Anderson, Hughes, Roos, and Seeger
\cite{AHRS2021} and Roos and Seeger \cite{RoosSeeger2023}.  The remaining
planar endpoint estimates were obtained by Lee, the author, F. Zhang, and Zhao
\cite{LRZZ2025}.

This restricted-radius theory has a local smoothing counterpart.  At
frequency $2^j$ the wave is essentially constant on time intervals of
length $2^{-j}$, so what matters is the number of such intervals needed to
cover $E$.  Beltran, Roos, Rutar, and Seeger formulate the resulting
fractal-time local smoothing problem in terms of the Assouad spectrum
\cite{BRRS2025}; the work \cite{LRZZ2025} proves related endpoint estimates
for the planar wave operator.  This is a natural off-diagonal extension of
the transfer mechanism in Section~\ref{sec:transfer}, but it leads to a
separate phase diagram and will not be pursued here.

\section{Open problems and concluding remarks}
\label{sec:open}

The preceding results leave a set of questions.  It is
helpful to separate the missing interior estimates, the two boundary
segments, and the geometric extensions of the problem.

\subsection{The remaining higher-dimensional interval}

Let $n\ge3$.  The only possible gap between the known sufficient and
necessary open regions is
$$
 2<p<P_n.
$$
The Euclidean local smoothing conjecture would close it.\footnote{Recall from the footnote in Section~\ref{sec:intro}
that proofs of this conjecture have been claimed in AI-generated
manuscripts, not yet verified by the community.}  For
$2<p<q_n$, the conjectural input is local smoothing with every positive
regularity $s>0$.  Substitution in Proposition~\ref{prop:transfer} gives
$$
 \operatorname{Re}\alpha>\frac1p-\frac{n-1}{2}.
$$
For $q_n\le p<P_n$, the conjectural regularity is
$s>(n-1)/2-n/p$, and the two terms in
\eqref{eq:transfer-condition} coincide, giving
$$
 \operatorname{Re}\alpha> -\frac{n-1}{p}.
$$
These are the two branches determined by necessity.  In particular, an
estimate at the Euclidean threshold $q_n$, followed by interpolation with the
$L^2$ identity below $q_n$ and with standard estimates above it, would
generate all the missing local smoothing inputs.

A potentially
weaker objective is to prove only the frequency-localized local smoothing
estimates with the special radial symbols produced by the Bessel expansion.
Such an estimate need not hold for every order-zero amplitude or every FIO
with cinematic curvature.  Determining the weakest wave estimate that still
implies the sharp maximal theorem would locate more precisely where the
difficulty of Stein's problem lies.

\subsection{The two negative-order boundary segments}

The equality cases should not be treated as a single endpoint problem.  On
the first segment,
$$
 2<p\le q_n,
 \qquad
 \operatorname{Re}\alpha=\frac1p-\frac{n-1}{2},
$$
the radial operator for real $\alpha$ is of restricted weak type but not of
weak type  \cite[Corollary~1.6]{NRS2023}.  Consequently neither weak nor strong type can hold for general
functions.  The remaining natural question on this segment is whether a
general restricted weak-type estimate is valid.  At the corner $p=q_n$ the
same distinction persists.

The second segment has a different character:
$$
 p>q_n,
 \qquad
 \operatorname{Re}\alpha=-\frac{n-1}{p}.
$$
Here the radial operator is already of strong type, so radial testing supplies
no obstruction to a general strong estimate.  The local smoothing transfer
does not reach equality because its arbitrarily small derivative loss becomes
a strict inequality in $\operatorname{Re}\alpha$.  Proving or disproving
strong type on this segment is therefore a maximal endpoint
question.  It may require estimates that provide quantitative control as the
local smoothing loss tends to zero, or an argument that does not pass through
Sobolev embedding in the dilation parameter.

For complex orders there is an additional issue of uniformity in
$\operatorname{Im}\alpha$.  Necessity sees only
$\operatorname{Re}\alpha$, whereas interpolation and the Bessel expansion
also produce constants depending on the imaginary part.  An endpoint
formulation should specify whether $\alpha$ is fixed or whether dependence on $|\operatorname{Im}\alpha|$ is required.

\subsection{Can the maximal operator be treated directly?}

The proof of Proposition~\ref{prop:transfer}  separates the
problem into the study of two half-waves separately, and it converts an $L^p$ norm in time
to a supremum by a one-dimensional Sobolev inequality.  A direct analysis of
the linearized maximal operator might instead provide the relation
between Bessel terms and the precise dependence of
their amplitudes on $t|\xi|$.

Several intermediate questions are therefore natural.  Can square-function
or variation estimates in $t$ replace the scalar Sobolev step?  Can a
time-frequency decomposition exploit that the phase is exactly
$t|\xi|$, rather than a general cinematic-curvature phase?  Can the maximal
operator be linearized in $t=t(x)$ and estimated by oscillatory-integral
methods that use the radial form of its symbol?  A positive answer could settle
Stein's problem without proving the full Euclidean local smoothing conjecture.
Conversely, any strategy showing that the sharp maximal estimate implies a
substantial local smoothing statement would establish that the present link
is deeper than a one-way transfer principle.

\subsection{Variable coefficients and manifolds}

Section~\ref{sec:dictionary} explains why the high-frequency pieces of
$M^\alpha$ belong to the cinematic-curvature FIO class.  The open geometric
problem is to construct the maximal family itself: Given a smooth family of curved hypersurfaces
$\Sigma_{x,t}$, one would like a holomorphic family of distributions
$m_{x,t}^\alpha$ such that $m_{x,t}^0$ is surface measure on
$\Sigma_{x,t}$ and the principal FIO has order
$-(n-1)/2-\operatorname{Re}\alpha$.  The construction must also control the
low-frequency contribution and all lower-order terms uniformly in $t$.

Once such a family is defined, the theorem of Gan, He, Li, and Wu suggests a
local maximal theorem for $t$ in a compact interval.  In odd dimensions their
threshold $P_n=p_n^+$ is already the sharp threshold for general
variable-coefficient local smoothing, in view of the examples of Minicozzi
and Sogge.  This does not by itself determine the maximal phase diagram:
necessity must be reformulated for the geometry of the chosen hypersurface
family, and endpoint behavior may depend on additional features of the
canonical relation.

The global problem is still more delicate.  Euclidean dilation invariance is
what passes from $[1,2]$ to all $t>0$ in
Proposition~\ref{prop:transfer}.  A manifold has no comparable global scaling,
and large times may encounter conjugate points, changes in topology, or
multiple branches of the geodesic flow.  Local and global generalized maximal
theorems should therefore be formulated separately.  The framework of
Mockenhaupt, Seeger, and Sogge and of Schlag and Sogge
\cite{MSS1992,SS1997} provides the natural starting point for the local
question.

\subsection{The classical endpoint as a guide}

The endpoint of the ordinary spherical maximal theorem illustrates why the
word ``boundary'' must be accompanied by a choice of function spaces.  Put
$\alpha=0$ and $p_0=n/(n-1)$.  If $n\ge3$, Bourgain proved
\cite{Bourgain1985} that the full spherical maximal operator satisfies
$$
 M^0:L^{p_0,1}(\R^n)\longrightarrow L^{p_0,\infty}(\R^n),
$$
although weak type $(p_0,p_0)$ fails.  In dimension two, even the general
restricted weak-type estimate at $p_0=2$ fails by a Besicovitch-set
construction, whereas it remains true for radial functions; see
\cite{RoosSeeger2026} for a recent account. For $n\ge3$ this endpoint lies below $p=2$ and is separate from the
negative-order boundary studied in most of this survey.  In the plane it is
the point $(p,\alpha)=(2,0)$ at the left end of the first boundary segment.

\subsection{Conclusion}

Stein's ``mystery'' has become a program.  The necessary
boundary is known, the radial theory is complete, the open planar region is
settled, and the new local smoothing theorems reduce the possible
higher-dimensional gap to $2<p<P_n$.  What remains distinguishes two directions:
an interior problem tied closely to Euclidean local smoothing, and equality (or endpoint)
cases for which local smoothing with an arbitrarily small loss is 
insufficient. Stein's analytic family continues to provide a meeting point for these questions because it records, in
a single complex parameter, the transition from averaging to differentiation
and from fixed-time curvature to space-time smoothing.

\section*{\texorpdfstring{Disclosure}{Disclosure}}

In preparing this manuscript, the author used ChatGPT Sol 6.1 (OpenAI) and
Claude Fable 5.1 (Anthropic) for language editing, for drafting and
reorganizing parts of the text from the author's notes, for checking
computations, and for locating relevant literature.  All mathematical
content has been verified by the author, who takes full responsibility
for it.

\medskip
\noindent{\textbf{Acknowledgments}.}
{The author is supported 
by the Spanish Agencia Estatal de Investigaci\'{o}n
through PID2023-146646NB-I00 funded by MICIU/AEI/10.13039/501100011033 and FEDER/EU and by ESF+ and through BCAM Severo Ochoa CEX2021-001142-S/MCIN/AEI/10.13039/501100011033 and by the Basque Government through the BERC 2026--2029 program.}

\end{document}